\documentclass[a4paper,12pt]{article}
\usepackage[english]{babel}
\usepackage[psamsfonts]{amssymb}
\usepackage{cmmib57}
\usepackage{exscale}
\usepackage{amsmath}
\usepackage{amscd}
\usepackage{latexsym}
\usepackage{amsthm}
\usepackage{amssymb}
\usepackage[latin1]{inputenc}
\usepackage[dvips]{graphicx}

\newcommand{\ott}{[0,T]}

\newcommand{\1}{{\bf 1}}

\newcommand{\R}{\mathbb R}

\newcommand{\cb}{\mathcal B}
\newcommand{\cac}{\mathcal C}

\newcommand{\cd}{\mathcal D}

\newcommand{\cf}{\mathcal F}
\newcommand{\ch}{\mathcal H}

\newcommand{\lc}{\left[}
\newcommand{\rc}{\right]}
\newcommand{\lcl}{\left\{}
\newcommand{\rcl}{\right\}}

\newcommand{\ga}{\gamma}

\newcommand{\si}{\sigma}

\newcommand{\al}{\alpha}
\newcommand{\Om}{\Omega}
\newcommand{\om}{\omega}
\newcommand{\del}{\delta}
\newcommand{\Del}{\Delta}

\newcommand{\Gam}{\Gamma}
\newcommand{\gam}{\gamma}
\newcommand{\ffi}{\varphi}
\newcommand{\hac}{\mathcal{H}}
\newcommand{\re}{\mathbb{R}}

\newcommand{\sig}{\sigma}
\newcommand{\ka}{\kappa}
\newcommand{\punt}{\cdot}
\newcommand{\tf}{\mathcal{F}}
\newcommand{\beq}{\begin{equation}}
\newcommand{\eeq}{\end{equation}}

\newtheorem{theorem}{Theorem}[section]

\newtheorem{definition}[theorem]{Definition}
\newtheorem{example}[theorem]{Example}

\newtheorem{lemma}[theorem]{Lemma}

\newtheorem{proposition}[theorem]{Proposition}

\newtheorem{remark}[theorem]{Remark}

\begin{document}

\title{{\bf The 1-d stochastic wave equation}\\ 
{\bf driven by a fractional Brownian motion}}
\author{
  \\
  {Llu\'{\i}s Quer-Sardanyons \thanks{Supported by an INRIA's Postdoc grant and grants BFM2003-01345 and HF2003-006, 
  Direcci\'on General de Investigaci\'on, Ministerio de Educaci\'on y Ciencia, Spain.}}              \\
  {\small\it Project OMEGA, INRIA Lorraine }  \\[-0.2cm]
  {\small\it Institut Elie Cartan, Campus Scientifique}          \\[-0.2cm]
  {\small\it BP 239 -- 54506 Vand\oe uvre-l\`es-Nancy, France}  \\[-0.2cm]
  {\small  {\tt quer@iecn.u-nancy.fr}}
  \\[-0.1cm]
\and
  { Samy Tindel}              \\
  {\small\it Institut Elie Cartan, Universit\'e de Nancy 1} \\[-0.2cm]
  {\small\it BP 239 -- 54506 Vand\oe uvre-l\`es-Nancy, France}  \\[-0.2cm]
  {\small  {\tt tindel@iecn.u-nancy.fr }} \\[-0.1cm]
  {\protect\makebox[5in]{\quad}}  
  \\
}

\maketitle

\begin{abstract}
In this paper, we develop a Young integration theory in dimension
2 which will allow us to solve a non-linear  
one dimensional wave equation
driven by an arbitrary signal whose rectangular increments satisfy
some H\"older regularity conditions, for some H\"older exponent
greater than $1/2$. This result will be applied to the infinite
dimensional fractional Brownian motion.
\end{abstract}

\vspace{2cm}


\vspace{0.3cm}


\newpage

\section{Introduction}

During the last past years, some spectacular advances have been made
in order to define and solve some stochastic PDEs driven by a Brownian
noise with a very general spatial covariance structure. This has been
achieved for instance in the infinite dimensional setting in
\cite{PZ,Pe} for both heat and wave equations, 
while the multiparametric setting has been treated first mainly for the 
1-d wave equation in \cite{walsh} and \cite{carmona}, and for the 2-dimensional case 
in \cite{DF} and \cite{MS}. The Walsh setting for 
SPDEs (\cite{walsh}) has been generalized then in \cite{Da}, leading to a quite complete 
picture for the stochastic heat equation, and also
for the wave equation up to dimension 3, allowing some deeper study 
in \cite{QS,qs2,DS}. Notice also
that an effort has been made in order to cover the case of L\'evy noises
e.g. in \cite{Sl,leon,erika2,mueller}.

\vspace{0.3cm}

On the other hand, since the recent introduction of fractional integrals
\cite{Za} and rough paths (see \cite{LQ,lejay}) techniques in
probability theory,
it has been clear that the pathwise method could be a good way to extend
the notion of stochastic differential equations to a wide variety
of noises. However, while this strategy has been explored thoroughly
in the case of ordinary differential equations (see e.g. \cite{CQ,NR}), the case of stochastic PDEs is still widely open. Indeed,
if the case of linear heat \cite{TTV} or wave \cite{Ca} 
equations driven by fractional noises
has been considered, sometimes leading to optimal results, 
only some very partial results are available in the case
of non-linear equations: let us mention for instance \cite{MN} for 
stochastic evolution equations driven by a fractional noise, with 
a very regular space covariance.
In this context, we have started to build  in \cite{GLT} an extension of 
Young integrals to the evolution setting, which allowed us to 
solve the stochastic heat equation driven by an infinite-dimensional
fractional Brownian noise with a non-trivial spatial covariance structure.
Observe that, while limited to the Young case, it seems that the methods
introduced in the latter reference can be extended to the rough path case
(see \cite{GT}).

\vspace{0.3cm}

The aim of the current paper is to show that the approach initiated in
\cite{GLT} can be extended to hyperbolic PDEs, and we will illustrate 
this fact by considering the stochastic wave equation in $\R$, of the form 
\beq\label{introwave}
\frac{\partial^2 Y}{\partial s^2}(s,t)
-\frac{\partial^2 Y}{\partial t^2}(s,t)=\sig(Y(s,t))\dot{X}(s,t),
\quad\mbox{ for }\quad
(s,t)\in [0,T]\times \re,
\eeq
with initial conditions given by
$$
Y(0,t)=\frac{\partial Y}{\partial s}(0,t)=0, 
\quad\mbox{ for }\quad
t\in \re.
$$
In equation (\ref{introwave}), $\sig$ stands for a smooth function
from $\re$ to $\re$, and $X$ represents the noise which drives the equation.
As usual in the SPDE theory, (\ref{introwave}) is understood
in the mild sense, which can be specified as follows: we will say that
$Y$ is a solution to (\ref{introwave}) if, for any $(s,t)\in\ott\times\re$,
we have
\beq
Y(s,t)=\int \int_{C(s,t)} \sig(Y(u,v)) X(du,dv),
\label{2.2i}
\eeq  
where $C(s,t)$ denotes the light cone with vertex $(s,t)$, that is the 
triangle delimited by the points $(s,t),(0,t+s)$ and $(0,t-s)$, and where 
the integral defining equation (\ref{2.2i}) is understood in the Young
sense. Then, for this latter equation,
we will give some existence and uniqueness results for a general 
class of noises whose rectangular increments
are H\"older continuous with H\"older exponent $> \frac12$
(see Theorem \ref{theorem1} for a precise statement).

\vspace{0.3cm}

One particular case of interest for us will be the  
infinite-dimensional fractional Brownian motion, which can be
defined in the following way: 
on a given complete probability space $(\Omega,\cf,P)$, let $X$ be
a centered Gaussian family 
$\{X(\phi);\phi\in\cd(\re^2)\}$
indexed by the set of test functions $\cd(\re^2)$, with a covariance 
function given by 
\beq\label{covfr}
E\left(X(\phi)X(\psi)\right)= 
c_H \int_{[0,T]^2} du dv |u-v|^{2H-2} 
\int_{\re^2}dx dy |x-y|^{-\nu} \phi(u,x)\psi(v,y),
\eeq
where $H>1/2$, $c_H=H(2H-1)$, and $\nu\in(0,1)$.
With a slight abuse of notation, set then 
$X(s,t)=X(\1_{[0,t]\times[0,s]})$.
Let also $\ch^{\gamma,\hat\gamma}(D)$
be the space of functions defined on $D\subset\R^2$ having a
H\"older regularity of order $\gamma$ in time and $\hat\gamma$ in space
(see Definition \ref{deffs} for the precise requirements). Eventually,
let $\mathcal{R}_{-\frac{\pi}{4}}$ be the -45° degree rotation in the plane.
In this context, our existence and uniqueness result
will be the following:
\begin{theorem}\label{fbmsol}
Let $X$ be a fractional Brownian noise defined by (\ref{covfr}). Suppose
that $\sig$ is a bounded function in $\cac^3(\re)$ with bounded derivatives of any order  
and let 
$\eta,\hat\eta\in(0,1)$
be such that $\eta+\hat\eta<1+H-\frac{\nu}{2}$. Then,
Equation (\ref{introwave}) admits a unique solution $Y$ such that 
$Y\circ \mathcal{R}_{-\frac{\pi}{4}}\in\ch^{\eta,\hat\eta}(R)$, for any rectangle $R$ around the origin.
\end{theorem}

\vspace{0.3cm}

As mentioned above, this paper can be seen as an extension
of \cite{GLT}, but the methodology used here
is quite different from the evolution type considerations 
contained in this latter reference. 
Indeed, it seemed easier, in the case of the wave equation, 
to consider the problem at hand in the multiparametric setting.
This has lead us to the following global strategy:
\begin{enumerate}
\item
Construct first a general Young integral on rectangles $R\subset\re^2$
whose sides are parallel to the axes. This integral will have the form
$\int_R f dg$, for two H\"older continuous functions $f,g:R\to\R$,
with large enough H\"older indexes. Notice that our construction is inspired
by \cite{Gu2}, but it is expressed here directly in terms of convergence of
Riemann sums, while \cite{Gu2} uses a 2-d analog of the $\Lambda$-map
defined in \cite{Gu}.
\item
Extend this Young integral in order to cover the case of a domain $R$
which is a triangle with two sides parallel to the axes. This is done
in a straightforward manner, by writing the triangle as a countable union of
rectangles.
\item
Rotate the wave equation in order to deal with an ordinary differential
equation in $\R^2$ involving the previous triangular domains. Once
our Young integral is constructed, the existence and uniqueness result will
be obtained by an extension of the usual fixed point argument for 
differential equations. It is worth noticing here that our computations
for this step will be quite delicate, in spite of having chosen a very 
regular coefficient $\si$. Indeed, though $\si\in\cac_b^3(\R)$, we 
will see that its interpretation as a map from $\ch^{\ga,\hat\ga}$
into itself does not enjoy the properties one usually assumes for the
resolution of Young equations: in fact, it is only locally Lipschitz
with quadratic growth, a fact which will add some technical difficulties to our analysis.
\item
In order to handle the case of the fractional Brownian noise, one has
to show that the rotation of this noise still satisfies the
H\"older regularity conditions allowing the definition of a Young
integral. This can be done in  our case, thanks to some almost explicit
and cumbersome calculations.
\end{enumerate}
This strategy will be made more explicit in the remainder of the
paper, but let us mention at this point that, to our
knowledge, Theorem \ref{fbmsol} is the first existence and 
uniqueness result for a non-linear wave equation driven by a general
kind of noise, and in particular by an infinite dimensional fractional
Brownian motion. We hope to extend this approach to a more irregular 
noise in a subsequent publication. Let us also mention that some of our
techniques can be related to those developped in \cite{Wa2} for numerical
approximation purposes.

\vspace{0.3cm}

Our paper will be structured as follows: at Section \ref{yint}, we will
define our general notion of Young integral in the plane. Then, we will 
solve the wave equation at Section \ref{weq}: Section \ref{extlight}
is devoted to the extension of the Young integral to the light cone.
We show how to rotate the wave equation at Section \ref{rotwa},
and then settle our fixed point argument at Section \ref{fixpo}.
The explicit application to the fractional Brownian  noise is left 
for Section \ref{applifbm}. Eventually, Section \ref{optim} has 
to be understood as a justification of the rotation trick for
our wave equation: we explore briefly another strategy consisting
in solving the equation, without previous rotation of the axes.
This leads to a great regularity loss of the Young integral, 
as well as some too restrictive 
assumptions on the driving noise.

\vspace{0.3cm}

Along the paper we will use the notation $C$ for any positive real constant, independently of its value.

\section{Two-dimensional Young integrals}\label{yint}

This section is devoted to a general result on Young integration in the
plane, which, to our knowledge, cannot be found in the literature,
in spite of being quite elementary:
we consider a rectangle $R=[s_1,s_2]\times [t_1,t_2]$,
where $s_1,s_2,t_1,t_2$ are arbitrary real numbers such that
$s_1<s_2$ and $t_1<t_2$,
 and we show that, 
under some regularity assumptions on the functions $x,y: R \rightarrow \re$, 
the integral $\int \int_R y(s,t) x(ds,dt)$ may be defined as a Young 
integral. 

\vspace{0.3cm}

Let us be more specific now about the regularity we will impose on 
the functions $x$ and $y$, and let us define the function spaces we
will consider in the sequel: first of all, for $\gam,\hat{\gam}\in (0,1)$,
set
\begin{equation}\label{defcgag}
\mathcal{C}^{\gam,\hat{\gam}}=\{f\in \mathcal{C}(\bar{R}), 
\|f\|_{\gam,\hat{\gam}}<\infty\},
\end{equation}
where
$$
\|f\|_{\gam,\hat{\gam}}= \sup_{s_1<s_2,t_1<t_2} 
\frac{|f(s_2,t_2)-f(s_2,t_1)-f(s_1,t_2)+f(s_1,t_1)|}
{|s_2-s_1|^\gam |t_2-t_1|^{\hat{\gam}}}.
$$  
We also suppose that the rectangle $R$ is contained in a sufficiently 
large square $\bar{R}=[R_1,R_2]^2$, which will be fixed 
throughout the discussion. With these notations in mind, the assumptions
on $x$ and $y$ will be the following:

\vspace{0.3cm}

\noindent {\bf{Hypothesis (H)}} The function $x$ belongs to the space 
$\mathcal{C}^{\gam,\hat{\gam}}$ 
and $y$ belongs to $\mathcal{C}^{\rho,\hat{\rho}}$, with 
$\gam+\rho>1$ and $\hat{\gam}+\hat{\rho}>1$. 
Moreover, there exist two positive constants $K, K'$ such that 
$$
|y(s,t)-y(s',t)|\leq K |s-s'|^\al, \; s,s'\in [s_1,s_2], \; t\in [t_1,t_2],$$
$$|y(s,t)-y(s,t')|\leq K' |t-t'|^\beta, \; s\in [s_1,s_2], \; t,t'\in [t_1,t_2],$$
with $\al>1-\gam$ and $\beta>1-\hat{\gam}$.

\vspace{0.3cm}

Let us also define the following functional spaces, in which the solutions
to our equations will live:
\begin{definition}\label{deffs}
For a function $y$ satisfying conditions (H), we define the semi-norm 
\begin{equation}\label{defno}
\|y\|:= \|y\|_{\rho,\hat{\rho}}+\|y\|_{1:\al}+\|y\|_{2:\beta},
\end{equation}
where the last two terms in the right-hand side denote the H\"older norms with respect to the first and second variable, respectively. Let then
$\mathcal{C}^{\rho,\hat{\rho}}_{\al,\beta}$ be the space of continuous 
functions $y$ such that 
$\|y\|<+\infty$, and observe that we will mostly consider the 
particular case $\hac^{\rho,\hat{\rho}}:= 
\mathcal{C}^{\rho,\hat{\rho}}_{\rho,\hat{\rho}}$.
\end{definition}

Let us describe now the discretization procedure we will use in order
to define our integral on $R=[s_1,s_2]\times [t_1,t_2]$:
for any rectangle $Q=[s,s']\times [t,t']$ and any function $g$ defined on $Q$, the rectangular increment of $g$ on $Q$ will be defined, as usual, by
$$\Del_Q g=g(s,t)-g(s,t')-g(s',t)+g(s',t').$$
For all $\del>0$, we consider $(\Pi^\del)_\del$ a family of partitions of the rectangle $R$ whose meshes goes to zero when $\del$ decreases to zero. Moreover, we assume that any of the partitions $\Pi^\del$ is formed by rectangles 
whose sides are parallel to the plane axes. Set $\Pi^\del=( (s^\del_i,t^\del_j))_{i,j}$, where $s_1=s^\del_0\leq s^\del_1\leq \dots \leq s^\del_{k^\del}=s_2$, $t_1=t^\del_0\leq t^\del_1\leq \dots \leq t^\del_{\bar{k}^\del}=t_2$.
With these notations in mind,
we consider the Riemann approximations
$$z^{\Pi^\del}_R=\sum_{i=0}^{k^\del-1}\sum_{j=0}^{\bar{k}^\del-1} y(s^\del_i,t^\del_j)\Del_{I^\del_{i,j}}x,$$
where we have used the notation $I^\del_{i,j}=[s^\del_i,s^\del_{i+1}]\times [t^\del_j,t^\del_{j+1}]$. 

\vspace{0.3cm}

Before stating our basic result on convergence
of Riemann sums, let us give an elementary property concerning the
partitions $\Pi^\del$:
\begin{lemma}
Let $R_1\leq s\leq t\leq R_2$ and let $s< r_1\leq\dots\leq r_k<t$ be a partition of $(s,t)$. Then, if $k\geq 2$, there exists an integer $l\in \{1,2,\dots,k\}$ such that
$$|r_{l+1}-r_{l-1}|\leq \frac{2}{k}|t-s|,$$
with the convention that $r_0=s$ and $r_{k+1}=t$.
\label{lemma}
\end{lemma}

\vspace{0.5cm}

\noindent
{\it{Proof}}. 
It is an immediate consequence of Lemma 2.2 in \cite{lejay}.
\hfill
$\Box$

\noindent

\vspace{0.3cm}

We are now in a position to state the main result of this section,
which gives the convergence of the Riemann sums defined above
to a limit $z_R=\int \int_{R} y(s,t) x(ds,dt)$:
\begin{proposition}
Recall that we have set $R=[s_1,s_2]\times [t_1,t_2]$. Then,
under Hypothesis (H), the sequence $\big(z^{\Pi^\del}_R\big)_\del$ converges, as $\del$ decreases to zero, to some limit denoted by $z_R$.
Furthermore, if we consider $z$ as a function of $s_1,s_2,t_1,t_2$, 
one gets that
\begin{equation}\label{remark1}
\left| \int \int_R y(s,t) x(ds,dt) \right| \leq C (\|y\|_\infty + \|y\|) \|x\|_{\gam,\hat{\gam}}  (s_2-s_1)^\gam (t_2-t_1)^{\hat{\gam}},
\end{equation}
and in particular, $z$
defines a continuous function
\begin{eqnarray*}
\bar{R}\times \bar{R} &\longrightarrow &\re\\
((s_1,t_1),(s_2,t_2))&\longmapsto& z_R.
\end{eqnarray*}
\label{proposition} 
\end{proposition}

\vspace{0.5cm}

\noindent
{\it{Proof }}. Fix $\del>0$ and $R=[s_1,s_2]\times [t_1,t_2]$. We will develop the proof in several steps, as follows.

\vspace{0.5cm}

\noindent
{\it{Step 1}}. We proceed, as in the proof of Proposition 2.1 in 
\cite{lejay}, by a kind of backward induction on the number of points
of the partition,
but instead of suppressing only one point, we will eliminate a whole column of $\Pi^\del$. Namely, owing to Lemma \ref{lemma}, we can choose an 
integer $\hat{\imath}\in \{1,2,\dots,k^\del-1\}$ such that
\beq
(s^\del_{\hat{\imath}+1}-s^\del_{\hat{\imath}-1})\leq \frac{2}{k^\del-1}(s_2-s_1).
\label{1}
\eeq
Consider now the new partition $\Pi$ of $R$ defined by
$$\Pi:
=\{(s^\del_i,t^\del_j), i=0,1,\dots,\hat{\imath}-1,\hat{\imath}+1,\dots,k^\del, j=0,1,\dots,\bar{k}^\del\}.
$$
Then, if we denote by $z^\Pi_R$ the Riemann sum corresponding to the partition $\Pi$, we obtain that
\begin{align}
z^{\Pi^\del}_R - z^\Pi_R  = & \sum_{j=0}^{\bar{k}^\del-1} y(s^\del_{\hat{\imath}-1},t^\del_j)\Del_{I^\del_{\hat{\imath}-1,j}}x + \sum_{j=0}^{\bar{k}^\del-1} y(s^\del_{\hat{\imath}},t^\del_j)\Del_{I^\del_{\hat{\imath},j}}x 
- \sum_{j=0}^{\bar{k}^\del-1} y(s^\del_{\hat{\imath}-1},t^\del_j)\Del_{[s^\del_{\hat{\imath}-1},s^\del_{\hat{\imath}+1}]
\times [t^\del_j,t^\del_{j+1}]}x \nonumber \\
 = &  \sum_{j=0}^{\bar{k}^\del-1} (y(s^\del_{\hat{\imath}},t^\del_j) - y(s^\del_{\hat{\imath}-1},t^\del_j)) \Del_{I^\del_{\hat{\imath},j}}x.
\label{2}
\end{align}
In order to get some upper bounds on $z^{\Pi^\del}_R - z^\Pi_R$, let us 
rewrite the last term in the above equality as a one-dimensional Riemann 
sum: set
$$h^\del_{\hat{\imath}} (t):=y(s^\del_{\hat{\imath}},t)-y(s^\del_{\hat{\imath}-1},t), \; t\in [t_1,t_2],$$
$$l^\del_{\hat{\imath}} (t):=x(s^\del_{\hat{\imath}+1},t)-x(s^\del_{\hat{\imath}},t), \; t\in [t_1,t_2].$$
Thus, with this notation, we get
$$z^{\Pi^\del}_R - z^\Pi_R = J^\del(\pi):= \sum_{j=0}^{\bar{k}^\del-1} h^\del_{\hat{\imath}} (t^\del_j)\left(l^\del_{\hat{\imath}} (t^\del_{j+1})- l^\del_{\hat{\imath}} (t^\del_j)\right),$$
where $\pi$ denotes the partition of $[t_1,t_2]$ given by  $\{t_1=t^\del_0\leq t^\del_1 \leq\dots \leq t^\del_{\bar{k}^\del}=t_2\}$. We will use now
the same kind of arguments as in \cite[Proposition 2.1]{lejay} 
to get suitable bounds on $J^\del(\pi)$. Indeed, applying Lemma \ref{lemma}, we can choose an integer $\hat{\jmath}\in \{1,\dots,\bar{k}^\del-1\}$ such that\beq
(t^\del_{\hat{\jmath}+1}-t^\del_{\hat{\jmath}-1})\leq \frac{2}{\bar{k}^\del-1}(t_2-t_1).
\label{3}
\eeq
One is then able to construct a new partition of $[t_1,t_2]$ 
in the following way: 
$$\tilde{\pi}:=\{t_1=t^\del_0\leq\dots\leq t^\del_{\hat{\jmath}-1}\leq t^\del_{\hat{\jmath}+1}\leq \dots\leq 
t^\del_{\bar{k}^\del}=t_2\}.$$
Hence, owing to the definition of $h$ and $l$, Hypothesis (H) and the bounds (\ref{1}) and (\ref{3}), we obtain
\begin{align*}
|J^\del(\pi)-J^\del(\tilde{\pi})| = & |h^\del_{\hat{\imath}}(t^\del_{\hat{\jmath}})-h^\del_{\hat{\imath}}(t^\del_{\hat{\jmath}-1})|\times
|l^\del_{\hat{\imath}}(t^\del_{\hat{\jmath}+1})-l^\del_{\hat{\imath}}(t^\del_{\hat{\jmath}})|\\
= & |\Del_{I^\del_{\hat{\imath}-1,\hat{\jmath}-1}}y|\times  |\Del_{I^\del_{\hat{\imath},\hat{\jmath}}}x|\\
\leq &  \|y\|_{\rho,\hat{\rho}} \|x\|_{\gam,\hat{\gam}}    (s^\del_{\hat{\imath}}-s^\del_{\hat{\imath}-1})^\rho (t^\del_{\hat{\jmath}}-t^\del_{\hat{\jmath}-1})^{\hat{\rho}}
 (s^\del_{\hat{\imath}+1}-s^\del_{\hat{\imath}})^\gam (t^\del_{\hat{\jmath}+1}-t^\del_{\hat{\jmath}})^{\hat{\gam}}\\
\leq &  C \frac{\|y\|_{\rho,\hat{\rho}} \|x\|_{\gam,\hat{\gam}}   }{(k^\del-1)^{\rho+\gam}(\bar{k}^\del-1)^{\hat{\rho}+\hat{\gam}}} (s_2-s_1)^{\rho+\gam} 
(t_2-t_1)^{\hat{\rho}+\hat{\gam}}.
\end{align*}
We can now proceed in a similar way to \cite[Proposition 2.1]{lejay}:
suppressing a carefully chosen point of $\tilde{\pi}$ 
and reiterating the process, one obtains, thanks to Hypothesis (H), that
$$|J^\del(\pi)-h^\del_{\hat{\imath}}(t_1)\left( l^\del_{\hat{\imath}}(t_2)-l^\del_{\hat{\imath}}(t_1)\right)|\leq 
C \frac{\|y\|_{\rho,\hat{\rho}} \|x\|_{\gam,\hat{\gam}}} {(k^\del-1)^{\rho+\gam}} (s_2-s_1)^{\rho+\gam} (t_2-t_1)^{\hat{\rho}+\hat{\gam}}.$$
Thus, going back to (\ref{2}), we get that 
$$|z^{\Pi^\del}_R-z^\Pi_R - E^\del_{\hat{\imath}}|\leq C \frac{\|y\|_{\rho,\hat{\rho}} \|x\|_{\gam,\hat{\gam}} }{(k^\del-1)^{\rho+\gam}} (s_2-s_1)^{\rho+\gam} (t_2-t_1)^{\hat{\rho}+\hat{\gam}},$$
where $E^\del_{\hat{\imath}}=h^\del_{\hat{\imath}}(t_1)\left( l^\del_{\hat{\imath}}(t_2)-l^\del_{\hat{\imath}}(t_1)\right)$. Now, one may reiterate the 
process and suppress a suitable column of $\Pi$, so that we end up with 
\beq
\left|z^{\Pi^\del}_R- \sum_{j=0}^{\bar{k}^\del-1} y(s_1,t^\del_j)\Del_{[s_1,s_2]\times [t^\del_j,t^\del_{j+1}]}x
\right| \leq 
C \|y\|_{\rho,\hat{\rho}} \|x\|_{\gam,\hat{\gam}} (s_2-s_1)^{\rho+\gam} (t_2-t_1)^{\hat{\rho}+\hat{\gam}} 
 + \sum_{r=1}^{k^\del-1}|E^\del_{\hat{\imath}_r}|.
\label{3.5}
\eeq
Let us estimate the last term of the right-hand side of the above inequality: 
for any $r=1,\dots,k^\del-1$, the integer 
$\hat{\imath}_r$ is an element of $\{1,\dots,k^\del-1\}$ such that
$(s^\del_{\hat{\imath}_r+1}-s^\del_{\hat{\imath}_r-1})\leq \frac{2}{k^\del-r}(s_2-s_1)$. Hence, by Hypothesis (H), one easily gets that
$$\sum_{r=1}^{k^\del-1}|E^\del_{\hat{\imath}_r}|\leq C \|y\|_{1:\al} \|x\|_{\gam,\hat{\gam}} (s_2-s_1)^{\gam+\al} (t_2-t_1)^{\hat{\gam}}.$$
Plugging this bound in (\ref{3.5}), we obtain that
\begin{align}
\left|z^{\Pi^\del}_R- \sum_{j=0}^{\bar{k}^\del-1} y(s_1,t^\del_j)\Del_{[s_1,s_2]\times [t^\del_j,t^\del_{j+1}]}x
\right| 
\leq  & C \|y\|_{\rho,\hat{\rho}} \|x\|_{\gam,\hat{\gam}} (s_2-s_1)^{\rho+\gam} (t_2-t_1)^{\hat{\rho}+\hat{\gam}} 
\nonumber \\
&  + C \|y\|_{1:\al} \|x\|_{\gam,\hat{\gam}}  (s_2-s_1)^{\gam+\al} (t_2-t_1)^{\hat{\gam}},
\label{4}
\end{align}
and we remark that the above bound (\ref{4}) is valid for any rectangle $R$ and any partition $\Pi^\del$ satisfying the underlying conditions. 

\vspace{0.5cm}

\noindent
{\it{Step 2}}.  Define now $Z^\del_{(s,t)}:=z^{\Pi^\del}_{R_0}$, where $R_0=[0,s]\times [0,t]$, for $s,t\in \re$ such that $R_0 \subset \bar{R}$. We aim to show that the sequence $(Z^\del_{(s,t)})_\del$ has a convergent subsequence. 
However, owing to Ascoli's Theorem, it is sufficient to prove that for any positive real number $K$, there exists $\eta>0$ such  that
\begin{equation}\label{asco}
\sup_{\|P-Q\|<\eta} |Z^\del_P-Z^\del_Q|\leq K,
\end{equation}
where $P=(s,t)$, $Q=(s',t')$.

\vspace{0.3cm}

In order to check the above condition, we consider $P=(s_1,t_1)$ and $Q=(s_2,t_2)$ two points satisfying $s_1<s_2$ and $t_1<t_2$, so that we can consider the rectangle $R=[s_1,s_2]\times [t_1,t_2]$. Let $(\Pi^\del)_\del$ be a family of partitions of the rectangle $[0,s_2]\times [0,t_2]$, formed by rectangles 
whose sides are parallel to the axes. In order to avoid tiresome notations, we denote again by $(s^\del_i,t^\del_j)$, $i=0,1,\dots k^\del$, $j=0,1,\dots, \bar{k}^\del$, the elements of $\Pi^\del$. Moreover, for the sake of simplicity we may assume that there exist two integers $1\leq n\leq k^\del$, $1\leq r\leq \bar{k}^\del$ such that $s^\del_n=s_1$ and $t^\del_r=t_1$. 
Then, owing to (\ref{4}) for the rectangles $R$, $[0,s^\del_{n-1}]\times[t^\del_{r+1},t_2]$ and  
$[s^\del_{n+1},s_2]\times[0,t^\del_{r-1}]$, we obtain that
\begin{align}
 & |Z^\del_Q-Z^\del_P|   =  \left|z^{\Pi^\del}_R + \sum_{i=0}^{n-2}\sum_{j=r+1}^{\bar{k}^\del-1} y(s^\del_i,t^\del_j) 
\Del_{I^\del_{i,j}}x + \sum_{i=n+1}^{k^\del-1}\sum_{j=0}^{r-2} y(s^\del_i,t^\del_j) 
\Del_{I^\del_{i,j}}x\right|\nonumber\\
 & \quad \leq  C \left\{ (s_2-s_1)^{\rho+\gam} (t_2-t_1)^{\hat{\rho}+\hat{\gam}} +  
(s_2-s_1)^{\gam+\al} (t_2-t_1)^{\hat{\gam}} \right\}
+ \left|\sum_{j=r}^{\bar{k}^\del-1} y(s_1,t^\del_j) 
\Del_{[s_1,s_2]\times [t^\del_j,t^\del_{j+1}]} x\right|    \nonumber\\
&   \quad \; +  C \left\{ (s^\del_{n-1})^{\rho+\gam} (t_2-t^\del_{r+1})^{\hat{\rho}+\hat{\gam}}+ 
(s^\del_{n-1})^{\gam+\al} (t_2-t^\del_{r+1})^{\hat{\gam}} \right\}
+ \left|\sum_{j=r+1}^{\bar{k}^\del-1} y(0,t^\del_j) 
\Del_{[0,s_1]\times [t^\del_j,t^\del_{j+1}]} x\right|\nonumber\\
&  \quad \; + C  \left\{ (s_2-s^\del_{n+1})^{\rho+\gam} (t^\del_{r-1})^{\hat{\rho}+\hat{\gam}} + 
(s_2-s^\del_{n+1})^{\gam+\al} (t^\del_{r-1})^{\hat{\gam}} \right\}
+ \left|\sum_{j=0}^{r-1} y(s_1,t^\del_j) 
\Del_{[s_1,s_2]\times [t^\del_j,t^\del_{j+1}]} x\right|.
\label{4.5}
\end{align} 
Let us bound now the terms in the right-hand side of (\ref{4.5}): first,
the sum of all the terms containing 
products of powers of $(s^\del_i-s^\del_{i'})$'s and 
$(t^\del_j-t^\del_{j'})$'s can be bounded, up to constants, by
\begin{align}\label{bnpd}
& (s_2-s_1)^{\rho+\gam} (t_2-t_1)^{\hat{\rho}+\hat{\gam}} + (s_2-s_1)^{\gam+\al} (t_2-t_1)^{\hat{\gam}}\\
& + (t_2-t_1)^{\hat{\rho}+\hat{\gam}} + (t_2-t_1)^{\hat{\gam}}
 + (s_2-s_1)^{\rho+\gam} + (s_2-s_1)^{\gam+\al} .\nonumber
\end{align}
Hence, we are left with the terms 
in the right-hand side of (\ref{4.5}) involving sums. 
Let us sketch the calculations for the first of these terms, namely
$$
S^\del\equiv 
\sum_{j=r}^{\bar{k}^\del-1} y(s_1,t^\del_j) 
\Del_{[s_1,s_2]\times [t^\del_j,t^\del_{j+1}]} x,
$$
since the two remaining terms can be treated analogously. 
On the other hand, $S^{\del}$ can be rewritten
as a Riemann sum, in the following way: 
$$S^\del=
\sum_{j=r}^{\bar{k}^\del-1} y(s_1,t^\del_j) (l(t^\del_{j+1})-l(t^\del_j)),
$$
where $l(t)=x(s_2,t)-x(s_1,t)$. Then we use the suppressing point argument, 
as it has been done at Step 1, and the regularity properties of $x$ and $y$. 
This easily yields 
\begin{align}\label{bnsd}
&|S^\del|\leq  C \|y\|_{2:\beta} \|x\|_{\gam,\hat{\gam}} (s_2-s_1)^\gam (t_2-t_1)^{\hat{\gam}+\beta}+ |y(s_1,t_1)||l(t_2)-l(t_1)|\\
& \leq C \|y\|_{2:\beta} \|x\|_{\gam,\hat{\gam}} (s_2-s_1)^\gam (t_2-t_1)^{\hat{\gam}+\beta}+ C \|y\|_\infty \|x\|_{\gam,\hat{\gam}}(s_2-s_1)^\gam (t_2-t_1)^{\hat{\gam}}.\nonumber
\end{align}
Hence, plugging (\ref{bnpd}) and (\ref{bnsd}) into (\ref{4.5}), we get
an upper bound of the form
$$
|Z^\del_P-Z^\del_Q|\le 
C \lc (s_2-s_1)^\mu + (t_2-t_1)^\xi \rc,
$$
with $\mu, \xi> 0$, from which inequality (\ref{asco}) easily follows.
Thus, owing to Ascoli's Theorem, there exists a subsequence of $(Z^\del)_\del$ converging uniformly to some continuous function $Z$ on $\bar{R}$. We make an abuse of notation and we denote also by $(Z^\del)_\del$ the underlying subsequence. 

\vspace{0.3cm}

Let us go back now to the definition of $z$, and for
$R=[s_1,s_2]\times[t_1,t_2]$,
set $z_R=Z_P+Z_Q-Z_{(s_1,t_2)}-Z_{(s_2,t_1)}$. Let us show that 
$z^{\Pi^\del}_R$ converges to $z_R$, when $\del$ tends to zero. 
For this, notice that we have the following decomposition:
$$
z^{\Pi^\del}_R=  z^{\Pi^\del}_Q-z^{\Pi^\del}_P-z^{\Pi^\del}_{[0,s_1]\times [t_1,t_2]}- 
z^{\Pi^\del}_{[s_1,s_2]\times [0,t_1]}.
$$
Thus, taking limit as $\del$ tends to zero in the above expression,  
we obtain that $z^{\Pi^\del}_R \rightarrow z_R$.
Furthermore, 
by continuity of $Z$, we deduce that the map $(P,Q)\mapsto z_R$ is continuous. 

\vspace{0.5cm}

\noindent
{\it{Step 3}}. 
Let us check that the limit of $( z^{\Pi^\del}_R)_\del$ is unique. 
This will be proved first
when considering a particular integrand, namely 
a function $\chi$ defined by some rectangular increment of 
a given path  $y$. Secondly, in the next Step 4, we will express $z^{\Pi^\del}_R$ in terms of a  
Riemann sum with respect to the function $\chi$ and other suitable terms, which will finally lead us to the uniqueness of the underlying sequence.

\vspace{0.3cm}

We will make use of the following fact: by a slight elaboration of the calculations done in the preceding Step 1 and Hypothesis (L), 
it can be proved that, for a rectangle $R=[s_1,s_2]\times[t_1,t_2]$,
we have
\begin{align}
 & |z^{\Pi^\del}_R -y(s_1,t_1)\Del_R x| \leq C  \|x\|_{\gam,\hat{\gam}} \Big(  \|y\|_{\rho,\hat{\rho}} 
 (s_2-s_1)^{\gam+\rho} (t_2-t_1)^{\hat{\gam}+\hat{\rho} } \nonumber \\ 
& \quad \quad +  \|y\|_{1:\al}  (s_2-s_1)^{\gam+\al} (t_2-t_1)^{\hat{\gam}}
 +  \|y\|_{1:\beta}  (s_2-s_1)^\gam (t_2-t_1)^{\hat{\gam}+\beta} \Big).
\label{s3}
\end{align}
For a given function $y$ satisfying Hypothesis (H), define then
the function $\chi:R\rightarrow \re$ as follows:
\beq
\chi(s,t):= \Del_{[s_1,s]\times [t_1,t]} y, \; (s,t)\in R,
\label{s3s3}
\eeq
and notice that $\chi$ depends on our particular choice of rectangle $R$.
Let us first study the regularity properties of $\chi$: for any rectangle 
$R_0$ contained in $R$, 
it can be easily checked that $\Del_{R_0} \chi =\Del_{R_0} y$. Owing to Hypothesis (H), this implies that $\chi \in \mathcal{C}^{\rho,\hat{\rho}}$.    
On the other hand, if $s, s'\in [s_1,s_2]$, $s<s'$, and $t\in [t_1,t_2]$, it 
holds that
$$|\chi(s,t)-\chi(s',t)|=|\Del_{[s,s']\times [t_1,t]} y|\leq C (s'-s)^\rho (t_2-t_1)^{\hat{\rho}}.$$
Thus, the function $\chi(\punt,t)$ is $\rho-$H\"older continuous uniformly 
with respect to $t$ and, moreover, one has that 
$\|\chi\|_{1:\rho} \leq C (t_2-t_1)^{\hat{\rho}}$. 
Analogously, it turns out that $\chi(s,\punt)$ is $\hat{\rho}-$H\"older 
continuous uniformly with respect to $s$ satisfying 
$\|\chi\|_{2:\hat{\rho}} \leq C (s_2-s_1)^\rho$.
Summing up, the function $\chi$ belongs to the space 
$\mathcal{H}^{\rho,\hat{\rho}}$. 

\vspace{0.3cm}

Hence, the calculations carried out in the preceding Steps 1 and 2 
hold true if we replace $y$ by $\chi$. In particular, by (\ref{s3}) and the above bounds for the H\"older norms of $\chi$, we have the following 
estimation:
\beq
\left| z^{\Pi^\del}_{\chi,R} \right| \leq K 
 (s_2-s_1)^{\gam+\rho} (t_2-t_1)^{\hat{\gam}+\hat{\rho}},
\label{s3bis} 
\eeq
where we have denoted by  $z^{\Pi^\del}_{\chi,R}$ the Riemann sum corresponding to the function $\chi$, that is 
$$z^{\Pi^\del}_{\chi,R}=\sum_{i=0}^{k^\del-1}\sum_{j=0}^{\bar{k}^\del-1} \chi(s^\del_i,t^\del_j)\Del_{I^\del_{i,j}}x.$$
Moreover, owing to Step 2, we obtain that the corresponding sequence $(Z^\del)_\del$ has a convergent subsequence and therefore deduce that 
$z^{\Pi^\del}_{\chi,R}$ converges, as $\del$ decreases to zero, to some limit $z_R$; notice that, in order to simplify notation, we 
do not point out the dependence of $\chi$ in  $Z^\del$ and $z_R$.

\vspace{0.3cm}

Let us check that the limit of $ (z^{\Pi^\del}_{\chi,R})_\del$ is unique. For this, we follow the same lines as in the proof of Proposition 2.1 in 
\cite{lejay}: let $\tilde{Z}$ be another limit of the sequence $(Z^\del)_\del$ and set $\tilde{z}_R := \tilde{Z}_P+\tilde{Z}_Q-\tilde{Z}_{(s_1,t_2)}-\tilde{Z}_{(s_2,t_1)}$ (recall that $P=(s_1,t_1)$ and $Q=(s_2,t_2)$). 
By (\ref{s3bis}), we obtain that
$$|z_R - \tilde{z}_R|\leq 2 K (s_2-s_1)^{\gam+\rho} (t_2-t_1)^{\hat{\gam}+\hat{\rho}},$$
which is indeed true for any rectangle $R$.
Thus, for any partition $\Pi_0=\{ R_{ij}=[s_i,s_{i+1}]\times [t_j,t_{j+1}], i=1,\dots,k, j=1,\dots,l\}$ of $R$, the following relation holds true: 
\begin{align*}
|z_R - \tilde{z}_R|& \leq \sum_{i=1}^k \sum_{j=1}^l |z_{R_{ij}}-\tilde{z}_{R_{ij}}| \\
& \leq 2 K \sum_{i=1}^k \sum_{j=1}^l (s_{i+1}-s_i)^{\gam+\rho} (t_{j+1}-t_j)^{\hat{\gam}+\hat{\rho}}\\
& \leq 2 K(s_2-s_1)(t_2-t_1) \left( \sup_{i} (s_{i+1}-s_i)^{\gam+\rho -1}\right) 
\left( \sup_{j} (t_{j+1}-t_j)^{ \hat{\gam}+\hat{\rho}-1}\right).
\end{align*}
Since $\gam + \rho>1$ and $\hat{\gam} + \hat{\rho}>1$, the above 
supremums tend to zero as the mesh of $\Pi_0$ decreases to zero,
which proves that $\tilde{z}_R=z_R$. 
Therefore, the limit of $(z^{\Pi^\del}_{\chi,R})_\del$ is unique.

\vspace{0.5cm}

\noindent
{\it{Step 4}}. Now we will show that, going back to the notations of Step 3, the sequence  $(z^{\Pi^\del}_R)_\del$ has a unique limit. Recall that  
$$z^{\Pi^\del}_R=\sum_{i=0}^{k^\del-1}\sum_{j=0}^{\bar{k}^\del-1} y(s^\del_i,t^\del_j)\Del_{I^\del_{i,j}}x,$$
where $x$ and $y$ satisfy Hypothesis (H). Then, the key point of our strategy
is to decompose $z^{\Pi^\del}_R$ in the following straightforward way:
\beq
z^{\Pi^\del}_R=  z^{\Pi^\del}_{\chi,R} + \sum_{j=0}^{\bar{k}^\del-1} y(s_1,t^\del_j)\Del_{[s_1,s_2]\times [t^\del_j,t^\del_{j+1}]}  x 
 + \sum_{i=0}^{k^\del-1} y(s^\del_i,t_1)\Del_{[s^\del_i,s^\del_{i+1}]
\times [t_1,t_2]}x - y(s_1,t_1)\Del_R x,
\label{dcp:zr}
\eeq
where the function $\chi$ is defined as in (\ref{s3s3}). 
Now, the uniqueness of the limit of $z^{\Pi^\del}_{\chi,R}$ has been
established in the previous step.
On the other hand, owing to Hypothesis (H), we are able to apply the convergence results in the one-dimensional setting (see \cite{young}) in order to 
obtain that 
$$\lim_{\del \searrow 0} \left(\sum_{j=0}^{\bar{k}^\del-1} y(s_1,t^\del_j)\Del_{[s_1,s_2]\times [t^\del_j,t^\del_{j+1}]}  x \right)=
\int_{t_1}^{t_2} y(s_1,v) d\left( x(s_2,v)-x(s_1,v)\right),$$
$$\lim_{\del \searrow 0} \left(\sum_{i=0}^{k^\del-1} y(s^\del_i,t_1)\Del_{[s^\del_i,s^\del_{i+1}]\times [t_1,t_2]}x \right)=
\int_{s_1}^{s_2} y(u,t_1) d\left( x(u,t_2)-x(u,t_1)\right),$$
where these limits are uniquely determined as one-dimensional Young
integrals. Going back to relation (\ref{dcp:zr}), this 
finishes the proof of the
uniqueness of the limit for the sequence $(z^{\Pi^\del}_R)_\del$.  
Moreover, the following relation is fulfilled:
\begin{align*}
\int \int_R y(u,v) dx(u,v)= & \int \int_R \left( \Del_{[s_1,u]\times [t_1,v]} y \right) dx(u,v)
 + \int_{t_1}^{t_2} y(s_1,v) d\left( x(s_2,v)-x(s_1,v)\right)\\
& + \int_{s_1}^{s_2} y(u,t_1) d\left( x(u,t_2)-x(u,t_1)\right) 
 - y(s_1,t_1)\Del_R x.
\end{align*}

\vspace{0.5cm}

\noindent
{\it{Step 5}}.
Eventually, owing to (\ref{s3}), it is readily checked that
$$
\left| \int \int_R y(s,t) x(ds,dt) \right| \leq C (\|y\|_\infty + \|y\|) \|x\|_{\gam,\hat{\gam}}  (s_2-s_1)^\gam (t_2-t_1)^{\hat{\gam}},
$$
for any rectangle $R=[s_1,s_2]\times [t_1,t_2]$, which ends the proof.
\qed 

\begin{remark}
Assume that the functions $x$ and $y$ satisfy the same hypothesis as in the previous Proposition \ref{proposition} and 
let $R=[s_1,s_2]\times [t_1,t_2]$ be a rectangle. Then,
as a concequence of Equation (\ref{s3}), we obtain the following estimation, which will be repeatedly applied throughout the proof of the existence and uniqueness of solutions to the wave equation (Theorem \ref{theorem1}):
\begin{align}
& \left| \int \int_R y(u,v) x(du,dv)\right|\leq C \|x\|_{\gam,\hat{\gam}} 
\Big\{ \|y\|_\infty (s_2-s_1)^\gam (t_2-t_1)^{\hat{\gam}} \nonumber \\
&\quad +  \|y\| \left( 
 (s_2-s_1)^{\gam+\rho} (t_2-t_1)^{\hat{\gam}+\hat{\rho}}  
+  (s_2-s_1)^{\gam+\al} (t_2-t_1)^{\hat{\gam}}
 +   (s_2-s_1)^\gam (t_2-t_1)^{\hat{\gam}+\beta} \right)\Big\}.
\label{afitacio}
\end{align}
\label{obs_afitacio}
\end{remark}

\section{The wave equation}\label{weq}

\begin{figure}[!ht]
\begin{center}
\scalebox{0.40}{\includegraphics{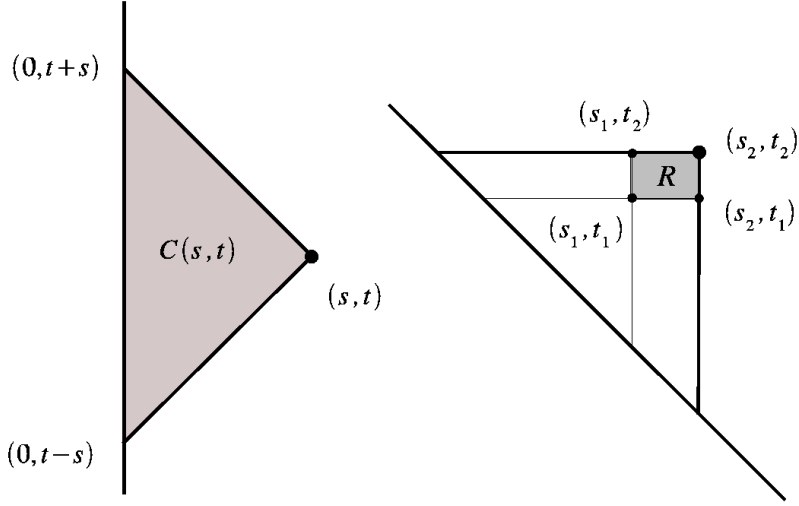}}   
\end{center}
\caption{On the left, the backward light cone with apex $(s,t)$ is represented, while on the right one can see a graphical representation of a 
rectangular increment of the rotated light cone $\tilde{C}$.}
\label{figure1}
\end{figure}

Let us turn now to the equation of main interest for us, that is
the following formal version of a perturbed wave equation
\beq\label{2.1}
\frac{\partial^2 Y}{\partial s^2}(s,t)
-\frac{\partial^2 Y}{\partial t^2}(s,t)=\sig(Y(s,t))\dot{X}(s,t),
\quad\mbox{ for }\quad
(s,t)\in [0,T]\times \re,
\eeq
with initial conditions given by
$$Y(0,t)=\frac{\partial Y}{\partial s}(0,t)=0, \; t\in \re.$$
Recall that we assume that
the real-valued functions $X,Y$ are defined on $[0,T]\times \re$, where $T$ is a fixed positive real number, and that
the coefficient $\sig$ is a real-valued smooth function (whose exact smoothness will be specified later on). 
Recall also that
we give a rigorous meaning to equation (\ref{2.1}) by means of its mild 
formulation, as follows: we will say that the continuous function $Y$ is 
a solution to (\ref{2.1}) if for any $(s,t)\in\ott\times\R$, it satisfies
the relation
\beq
Y(s,t)=\int \int_{C(s,t)} \sig(Y(u,v)) X(du,dv),
\label{2.2}
\eeq  
where $C(s,t)$ denotes the open light cone with vertex $(s,t)$ and projected to the $t-$axle, that is the 
triangular domain delimited by the points $(s,t),(0,t+s)$ and $(0,t-s)$ (see Figure \ref{figure1}).
We will also assume that 
the integral defining Equation (\ref{2.2}) is understood in the Young
sense given by our Proposition \ref{proposition}.

\vspace{0.3cm}

Before going into the details of the definition of our equation,
let us specify first our assumptions on the function $\sig$.
In fact, the main property we will need on this coefficient
can be summarised as follows:

\vspace{0.3cm}

\noindent {\bf{Hypothesis (L)}} $\sig: \re\rightarrow \re$ is a smooth function preserving the regularity properties on spaces of the form 
$\hac^{\rho,\hat{\rho}}$, for $\rho,\hat{\rho}>0$. 
Moreover, it satisfies the following two conditions: 
\beq
\|\sig(y)\| \leq C \|y\| (1+\|y\|),
\label{cresig}
\eeq
\begin{align}
\|\sig(y_1)-\sig(y_2)\| \leq K & \left( \|y_1-y_2\|_\infty + \|y_1-y_2\|\right) \nonumber \\ 
& \; \times \left( 1+ \|y_1\|  + \|y_2\|+ \|y_1-y_2\| +  \left( \|y_1\|+ \|y_1-y_2\|\right)^2 \right),
\label{lipsig} 
\end{align}
for any $y,y_1,y_2 \in \hac^{\rho,\hat{\rho}}$ and some positive constants $K, C$.

\vspace{0.3cm}

\noindent
This assumption will be made throughout the paper, and one should 
observe that it is satisfied in the following simple case:
\begin{lemma}
Assume that $\sig$ is bounded, belongs to the space $\mathcal{C}^3 (\re)$ and has bounded derivatives. Then, Hypothesis (L) is fulfilled.
\label{lemma1}
\end{lemma}

\vspace{0.5cm}

\noindent
{\it{Proof}}. 
In this proof, we will use the same kind of arguments
as in \cite{GT}, and thus only the main ideas of our strategy
will be sketched. To begin with, 
we show that $\sig$ preserves the regularity on spaces of the form 
$\hac^{\rho,\hat{\rho}}$, which amounts to control all the norms used
in (\ref{defno}) to define $\|\sig(y)\|$. 
\vspace{0.3cm}

First of all, let us deal with the regularity of $\sig(y)$ on rectangles. 
For any rectangle $R=[s_1,s_2]\times [t_1,t_2]$, the following equality is fulfilled:
$$\Del_R \sig(y)=\int_0^1 dr \int_0^1 d\tau \partial_r \partial_\tau \sig(a(r,\tau)),$$
where
\beq
a(r,\tau)=y(s_1,t_1)+r(y(s_2,t_1)-y(s_1,t_1))+\tau(y(s_1,t_2)-y(s_1,t_1))+ r\tau\Del_R y.
\label{3.4}
\eeq
Thus,
$$\Del_R \sig(y)=\int_0^1 dr \int_0^1 d\tau \left[ \sig'(a(r,\tau)) \partial_r \partial_\tau a 
+\sig''(a(r,\tau)) \partial_r a \partial_\tau a\right],$$
and, as in \cite{GT}, one can deduce that
$$|\Del_R \sig(y)| \leq C  \|y\| (1+\|y\|)   (s_2-s_1)^\rho (t_2-t_1)^{\hat{\rho}}.$$
The $(1:\rho)$ and $(2:\hat{\rho})-$Hölder regularities follow from the Lipschitz property of $\sig$. Indeed, it is straitghforward to check that
$\|\sig(y)\|_{1:\rho}\leq C \|y\|_{1:\rho}$ and $\|\sig(y)\|_{2:\hat{\rho}}\leq C \|y\|_{2:\hat{\rho}}$.
Hence, we conclude that $\sig(y)$ belongs to $\hac^{\rho,\hat{\rho}}$ and, moreover, that condition (\ref{cresig}) is fulfilled.

\vspace{0.3cm}

On the other hand, we have to check that $\sig$ satisfies the local Lipschitz 
property (\ref{lipsig}).
Let us sketch the calculations concerning the $(1:\rho)-$Hölder norm; the $(2:\hat{\rho})-$Hölder norm 
may be carried out using the same arguments:
let $s,\bar{s},t$ be such that $(s,t), (\bar{s},t)$ belong to $\bar{R}$. Then, following the same lines as 
in \cite{GT}, it can be proved that
$$ \sig(y_1(s,t))-\sig(y_2(s,t))- \sig(y_1(\bar{s},t))+\sig(y_2(\bar{s},t))=
\int_0^1 dr \int_0^1 d\tau \partial_r \partial_\tau \sig(b^r (\tau)),$$
where $b^r(\tau)=y^r(\bar{s})+\tau(y^r(s)-y^r(\bar{s}))$ and $y^r=y_2+r(y_1-y_2)$. Expanding the right-hand side of the above expression and using the assumptions on $\sig$ and $y_1,y_2$, one ends up with
\beq
\| \sig(y_1)-\sig(y_2)\|_{1:\rho} \leq   C \|y_2\| \left( \| y_1-y_2\|_\infty + \| y_1-y_2\|\right)
 +  \| y_1-y_2\| \left( 1+ \| y_1-y_2\|_\infty\right).
\label{3.6}
\eeq
As it has been mentioned above, we have an analogous bound for the $(2:\hat{\rho})-$Hölder norm. 

Eventually, in order to deal with the regularity on rectangles
of $\sig(y_1)-\sig(y_2)$, 
we notice that the following equality holds true:
$$\Del_R (\sig(y_1)-\sig(y_2))=\int_0^1 dr \int_0^1 d\tau \int_0^1 d\nu \partial_r \partial_\tau
\partial_\nu \sig(a^r(\tau,\nu)),$$
with $a^r=a_1 +r(a_2-a_1)$ and $a_i$ defined as in (\ref{3.4}) but with $y$ replaced by $y_i$, $i=1,2$ (see \cite{GT}). 
Then, it can be shown that
\begin{align}
& \| \sig(y_1)-\sig(y_2)\|_{\rho,\hat{\rho}} \leq   C  \| y_1-y_2\| \nonumber \\
& \quad \quad + C \left( \| y_1-y_2\|_\infty + \| y_1-y_2\|\right) \left( \|y_1\| + \| y_1-y_2\|\right) \left( 1+ \|y_1\| + \| y_1-y_2\|\right).
\label{3.8}
\end{align}
Putting together the bounds (\ref{3.6}) and (\ref{3.8}), we conclude the proof. 
\qed

\vspace{0.3cm}

We are now ready to define rigorously our wave equation (\ref{2.2}).

\subsection{Extension of the integral to the light cone}\label{extlight}

A first step towards a rigorous definition of Equation (\ref{2.2})
is to extend slightly our definition of Young integral in order to cover the
case of a triangular domain like $C(s,t)$.
This will be done by a straightforward limiting argument,
as follows. 

\vspace{0.3cm}

Let $\mathcal{R}(s,t)$ be the set of families of rectangles $(R_n)_n$ 
of the form $[s^n_1,s^n_2]\times [t^n_1,t^n_2]$,
such that $\uplus_{n=1}^\infty R_n =C(s,t)$ and
\beq
\sum_{n=1}^\infty (s^n_2-s^n_1)^\gam  (t^n_2-t^n_1)^{\hat{\gam}}< +\infty.
\label{2.25}
\eeq
Then, if $X,Y$ satisfy Hypothesis (H) from Section 1, the integral 
$$\int \int_{R_n} \sig(Y(s,t)) X(ds,dt)$$
is well defined, for all $n\geq 1$. Moreover, by 
relations  (\ref{remark1}) and (\ref{2.25}),  for $(R_n)_n$ in $\mathcal{R}(s,t)$, the series 
\begin{equation}\label{defrec}
\int \int_{C(s,t)} \sig(Y(s,t)) X(ds,dt)
:=\sum_{n=1}^\infty \int \int_{R_n} \sig(Y(s,t)) X(ds,dt)
\end{equation}
is finite and the limit does not depend on the chosen element of $\mathcal{R}(s,t)$. From now on, the integral $\int \int_{C(s,t)} \sig(Y(s,t)) X(ds,dt)$
will be understood by means of (\ref{defrec}), which gives a reasonable
definition of a Young integral on $C(s,t)$.

\subsection{Rotation of the wave equation}\label{rotwa} 

As it will be pointed out at Section \ref{optim}, it is convenient
to deal with the existence and uniqueness of solution to Equation (\ref{2.1})
after a change of coordinates corresponding to a 45° rotation, which we proceed to detail now. 
Set $y(s,t):=Y\left(\frac{t+s}{\sqrt{2}},\frac{t-s}{\sqrt{2}}\right)$ and 
$x(s,t):=X\left(\frac{t+s}{\sqrt{2}},\frac{t-s}{\sqrt{2}}\right)$. Then, 
a trivial change of variables in 
the integral Equation (\ref{2.2}) yields that $y$ satisfies
\beq
y(s,t)=\int \int_{\tilde{C}(s,t)} \sig(y(u,v)) x(du,dv),
\label{2.3}
\eeq
where now $\tilde{C}(s,t)$ corresponds to the light cone with vertex $(s,t)$ and projected to the line $\{(s,-s), s\in \re\}$, that is the triangular
domain delimited by $(s,t)$, $(s,-s)$ and $(-t,t)$. 
Notice that the domain of definition of our original equation was included
in an arbitrary large rectangle $\bar{R}=[R_1,R_2]^2$, and let $\tilde{R}$
be the image of $\bar{R}$ under the rotation. Then  we
assume that the norm $\|x\|$ of $x$ can be controlled suitably on $\tilde{R}$,
and the domain of definition of the rotated equation will be 
\begin{equation}\label{domrot}
D_T:=\{(s,t)\in \tilde{R}, -s\leq t\leq -s+\sqrt{2}T\},
\end{equation} 
for a given arbitrary $T>0$. 
The new initial conditions are given by
\beq
y(s,-s)=0, \frac{\partial y}{\partial s}(s,-s)=- \frac{\partial y}{\partial t}(s,-s), \; s\in \re.
\label{2.4}
\eeq
Obviously, Equation (\ref{2.3}) assumes implicitely that $x$ and $y$
satisfy Hypothesis (H), and all our statements will make use of this 
hypothesis. Then we will show at Section \ref{applifbm} that this assumption
can be made when $X$ is an infinite dimensional fractional Brownian
motion.
\begin{remark}
Suppose that the functions $x,y$ satisfy Hypothesis (H) and that $\sig$ satisfies Hypothesis (L). Then, owing to Proposition \ref{proposition} 
and the considerations in the preceding Subsection \ref{extlight}, the integral of $y$ with respect to $x$, namely
$$I(s,t):=\int \int_{\tilde{C}(s,t)} \sig(y(u,v)) x(du,dv),$$
is well defined. Moreover, it can be easily seen that, for any rectangle 
$R=[s_1,s_2]\times [t_1,t_2]$, we have 
$$\Del_R I=\int \int_R \sig(y(u,v)) x(du,dv);$$
a graphical proof of this fact is given  in Figure \ref{figure1}.
Taking into account relation (\ref{remark1}), this implies that 
$I\in \ch^{\gam,\hat{\gam}}$. 
\label{remark2}
\end{remark}

\begin{example}
Let $(s,t)\in \tilde{R}$ and consider $\tilde{C}(s,t)$ the {\emph{rotated}} light cone with vertex $(s,t)$. Assume that $x,y$ satisfy Hypothesis (H) with $\gam +\hat{\gam}> 1$ (this condition will be implied by the assumptions in 
Theorem \ref{theorem1}). Assume also that $\sig$ satisfies Hypothesis (L). Then, the integral $I(s,t)$ is constructed as in Section 3.1. 
Indeed, it can be easily shown that 
$\tilde{C}(s,t)$ can be recovered by the family of rectangles 
$\uplus_{k=1}^\infty M_k$, where $M_k$ can  in turn be written as
a union of $2^{k-1}$ squares of side $\frac{t+s}{2^k}$.
Then, in this case we have that
\begin{align*}
\sum_{n=1}^\infty (s^n_2-s^n_1)^\gam  (t^n_2-t^n_1)^{\hat{\gam}} & =
\sum_{k=1}^\infty 2^{k-1} \left( \frac{t+s}{2^k}\right)^{\gam + \hat{\gam}}\\
& = \frac{(t+s)^{\gam + \hat{\gam}}}{2}   \sum_{k=1}^\infty \frac{1}{2^{k(\gam + \hat{\gam}-1)}}\leq C  
(t+s)^{\gam + \hat{\gam}},
\end{align*}
and we obtain the following estimate, which will be useful
in the sequel:
$$\left| \int \int_{\tilde{C}(s,t)} \sig(y(u,v)) x(du,dv) \right| \leq C \|x\|_{\gam,\hat{\gam}} \left(1+ \|y\|(1+ \|y\|)\right) (t+s)^{\gam +\hat{\gam}}.$$
\label{example1}
\end{example}

\subsection{Existence and uniqueness of solution}\label{fixpo}

We are now ready to state and prove our main general result on 
existence and uniqueness of the solution for Equation (\ref{2.2}), under the
form (\ref{2.3}). Recall that $ \tilde{R}$ denotes the image of 
$\bar{R}=[R_1,R_2]^2$ 
under the 45° rotation, and that $D_T$ has been defined by (\ref{domrot}).
\begin{theorem}
Assume that the function $x$ belongs to $\mathcal{C}^{\gam,\hat{\gam}}(D_T)$, for some $\gam,\hat{\gam} \in (\frac{1}{2},1)$, 
and that $\sig$ satisfies Hypothesis (L). Then, there exists a unique solution  to Equation (\ref{2.3}) in $\mathcal{H}^{\ka,\hat{\ka}}(D_T)$, 
for all $\ka \in (1-\gam,\gam)$ and $\hat{\ka}\in (1-\hat{\gam},\hat{\gam})$,  denoted by 
$\{y(s,t), (s,t)\in D_T\}$, satisfying the initial conditions (\ref{2.4}).  
\label{theorem1}
\end{theorem}

\begin{remark}
The previous theorem ensures that there exists a unique solution to Equation (\ref{2.3}) in the bounded domain 
$D_T=\{(s,t)\in \tilde{R}, -s\leq t\leq -s+\sqrt{2}T\}$. 
However, we aim to have a solution in the whole domain, namely  
$$\{y(s,t), (s,t)\in \mathbb{R}^2, -s\leq t\leq -s+\sqrt{2}T\}.$$
One may construct this solution as follows: let $(\tilde{R}_n)_n$ be a 
family of squares in $\mathbb{R}^2$ such that
\begin{itemize}
\item[(i)] $\tilde{R}_n \subset \tilde{R}_{n+1}$, for all $n$, 
\item[(ii)] $\cup_n \tilde{R}_n =\mathbb{R}^2$.
\end{itemize}  
We denote by $y_n$ the unique solution of Equation (\ref{2.3}) on $\{(s,t)\in \tilde{R}_n, -s\leq t\leq -s+\sqrt{2}T\}$. Then, it is readily checked that 
the following is fulfilled: if $n>m$,
$$y_n (s,t)=y_m (s,t), \;\textrm{for}\;  (s,t)\in \tilde{R}_m \;\;\textrm{satisfying}\;\; -s\leq t\leq -s+\sqrt{2}T.$$ 
This let us define our global solution without ambiguity: 
$y(s,t):= y_n(s,t)$, for all $(s,t)\in \tilde{R}_n$ with  $-s\leq t\leq -s+\sqrt{2}T$.
\end{remark}

\vspace{0.5cm}

\begin{figure}[!ht]
\begin{center}
\scalebox{0.40}{\includegraphics{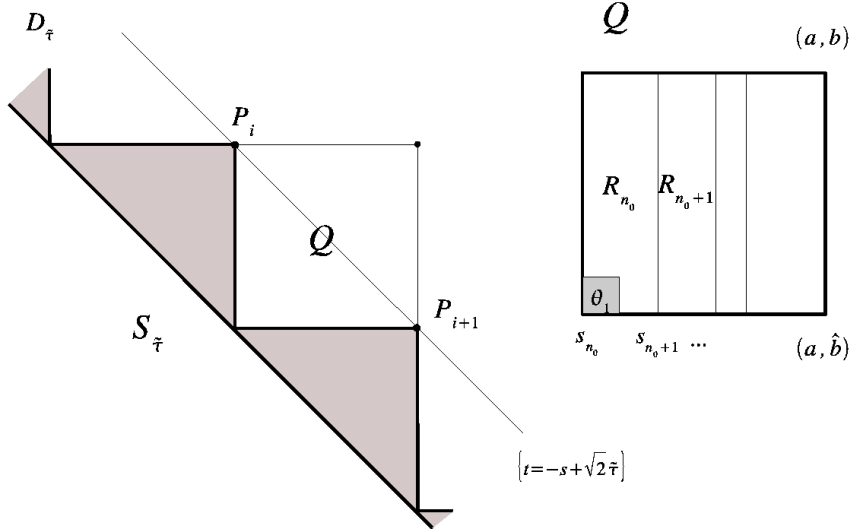}}   
\end{center}
\caption{On the left-hand side, a graphical description of the {\it{stairs}} domain $S_{\tilde{\tau}}$ is given, while on the right-hand side the decomposition of the square $Q$ in terms of the rectangles' sequence $(R_n)_n$ is represented.}
\label{figure2}
\end{figure}

\noindent
{\it{Proof of Theorem \ref{theorem1}}}.  Let us fix $\gam,\hat{\gam}>\frac{1}{2}$ and $\ka \in (1-\gam,\gam)$, $\hat{\ka}\in (1-\hat{\gam},\hat{\gam})$. 
We will make use of a fixed-point argument. For this, we define the map $\Gam$, from $\hac^{\ka,\hat{\ka}} (D_T)$
into itself, as follows:
\beq
\Gam (y)(s,t):= \int \int_{\tilde{C}(s,t)} \sig(y(u,v)) x(du,dv).
\label{2.5}
\eeq
Notice that $\Gam$ is well defined and, thanks to Remark \ref{remark2} and the fact that $\ka<\gam$ and $\hat{\ka}<\hat{\gam}$, $\Gam(y)\in \hac^{\ka,\hat{\ka}}(D_T)$. 
The strategy in order to 
show that $\Gam$ has a unique fixed point is quite standard: first we prove that $\Gam$ maps some closed ball of 
$\hac^{\ka,\hat{\ka}}(D_T)$ into itself and secondly, that it is a contraction. For the latter to be fulfilled, it is sufficient to show that
$$\| \Gam(y_1)-\Gam(y_2)\| \leq K \|y_1-y_2\|,$$
for all $y_1,y_2\in \mathcal{H}^{\ka,\hat{\ka}}(D_T)$ and some positive constant $K<1$.
We will develop the proof in several steps, which may be summarised as follows:
\begin{enumerate}
\item We will prove first the existence and uniqueness of a fixed point in a small part of the domain, 
namely in a band $D_{\tilde{\tau}}$, for some sufficiently small $\tilde{\tau}<T$. 
\item In order to iterate the procedure and cover the whole domain $D_T$, we will consider a {\it{stairs}} domain $S_{\tilde{\tau}}$ contained in $D_{\tilde{\tau}}$. 
Then we will proceed to prove the existence and uniqueness of a fixed point in one of the squares determined by two steps of the stairs, denoted by $Q$. For this, 
we will recover $Q$ by a suitable finite family of rectangles $(R_n)_n$ and prove that there exists a ball in $\mathcal{H}^{\ka,\hat{\ka}}(R_n)$
which is left invariant by $\Gamma$, for all $n$. In Figure \ref{figure2}, the {\it{stairs}} domain $S_{\tilde{\tau}}$, together with 
the square $Q$ and the sequence $(R_n)_n$, are represented.
This step will contain most of the technical difficulties of our proof,
since Hypothesis (L) only assumes that $\si$ is a locally Lipschitz
coefficient, which is usually considered as a too mild assumption
in the Young integration theory.
\item In Step 3 we proceed to show that the corresponding map defined on the invariant ball of $\mathcal{H}^{\ka,\hat{\ka}}(R_n)$ is a contraction. This will be carried out again by recovering $R_n$ 
by a convenient family of equally sized squares.
\item Finally, we will iterate this procedure in order to get a unique fixed point of $\Gam$ in a larger stairs domain than 
$S_{\tilde{\tau}}$. This will let us cover the whole band $D_T$.
\end{enumerate}

\vspace{0.5cm}

\noindent {\it{Step 1: Fixed point in a small part of the domain: the band $D_{\tilde{\tau}}$}}\\
We show first that there exists $\tau<T$ such that the closed ball of 
$\ch^{\ka,\hat{\ka}}(D_T)$, namely
$$
\mathcal{B}_{\tau,\bar{K}}
:=\{ y\in \hac^{\ka,\hat{\ka}}(D_\tau);\, 
y(s,-s)=0, \,\|y\|\leq \bar{K}\},
$$
is invariant under $\Gam$, for some $\bar{K}>0$.  For this,
recall that, under the standing assumptions,
$$\|y\|= \|y\|_{\ka,\hat{\ka}} + \|y\|_{1:\ka} + \|y\|_{2:\hat{\ka}},$$
and we will try to bound all the norms above separately:

\vspace{0.3cm}

First, let us deal with the Hölder norm in rectangles. If $R=[s_1,s_2]\times [t_1,t_2]$ is a rectangle included in $D_\tau$, then we know that
$$|\Del_R \Gam(y)|=\left| \int \int_R \sig(y(u,v)) x(du,dv)\right|.$$
Furthermore, owing to (\ref{s3}) and (\ref{cresig}), we have that,
\begin{align*}
& |\Del_R \Gam(y)-\sig(y(s_1,t_1))\Del_R x| \\ 
& \quad \leq C \|x\|_{\gam,\hat{\gam}} \|y\| (1+\|y\|) \Big( (s_2-s_1)^{\gam-\ka} (t_2-t_1)^{\hat{\gam}-\hat{\ka}}\Big)  
(s_2-s_1)^\ka (t_2-t_1)^{\hat{\ka}}\\
& \quad \leq C \|x\|_{\gam,\hat{\gam}} \|y\| (1+\|y\|) g(\tau) (s_2-s_1)^\ka (t_2-t_1)^{\hat{\ka}},
\end{align*}
where $g(\tau)$ tends to zero, as $\tau$ decreases to zero. Taking into account that 
$$
|\sig(y(s_1,t_1))\Del_R x|\leq 
\|\sig\|_\infty \|x\|_{\gam,\hat{\gam}} (s_2-s_1)^\gam (t_2-t_1)^{\hat{\gam}}
$$ 
and $\ka<\gam$, $\hat{\ka}<\hat{\gam}$, we obtain 
\beq
\|\Gam(y)\|_{\ka,\hat{\ka}} \leq C  g(\tau) (1+\|y\|( 1+\|y\|) ),
\label{2.65}
\eeq
where we still denote by $g(\tau)$ a positive function decreasing to zero, as $\tau$ tends to zero.  


\vspace{0.3cm}

Concerning the Hölder norm with respect to the first variable, we observe that, for any $s'<s$ and $t$ such that 
$(s',t),(s,t)\in D_\tau$, we have 
\begin{align}
|\Gam(y)(s,t)-\Gam(y)(s',t)|  \leq  & \left|\int \int_{[s',s]\times [-s',t]} \sig(y(u,v)) x(du,dv)\right| \nonumber \\ 
 & + \left|\int \int_{\tilde{C}(s,-s')} \sig(y(u,v)) x(du,dv)\right|.
\label{2.75}
\end{align}
On one hand, by the same calculations carried out to obtain (\ref{2.65}) or, equivalently, by (\ref{afitacio}), one easily gets that
\beq
\left|\int \int_{[s',s]\times [-s',t]} \sig(y(u,v)) x(du,dv)\right|\leq C  (1+\|y\| (1+\|y\|)) \tilde{g}(\tau) (s-s')^\ka,
\label{2.8}
\eeq
with $\tilde{g}$ converging to zero as $\tau\searrow 0$. On the other hand, 
owing to Example \ref{example1}, it holds that
\begin{align}
\left|\int \int_{\tilde{C}(s,-s')} \sig(y(u,v)) x(du,dv)\right|&\leq C \|x\|_{\gam,\hat{\gam}} (1+\|y\| (1+\|y\|))  (s-s')^{\gam+\hat{\gam}}\nonumber \\ 
& \leq C  (1+\|y\| (1+\|y\|)) (\sqrt{2} \tau)^{\gam-\ka+\hat{\gam}} (s-s')^\ka,
\label{2.9}
\end{align}
and plugging (\ref{2.8}) and (\ref{2.9}) in (\ref{2.75}), 
we obtain the following estimate:
\beq
\|\Gam(y)\|_{1:\ka} \leq C  (1+\|y\|(1+\|y\|) ) \tilde{g}(\tau),
\label{3.1}
\eeq
where we use again the same notation $\tilde{g}(\tau)$ for a function converging to zero as $\tau\searrow 0$. 

\vspace{0.3cm}

Using the same arguments as for the $(1:\ka)-$norm, one can also get the bound
\beq
\|\Gam(y)\|_{2:\hat{\ka}} \leq C  (1+\|y\| (1+\|y\|)) \bar{g}(\tau),
\label{3.2}
\eeq
with $\bar{g}$ satisfying again $\lim_{\tau\searrow 0}\bar{g}(\tau)=0$. Therefore, 
putting to\-ge\-ther 
the three bounds (\ref{2.65}), (\ref{3.1}) and (\ref{3.2}), we end up with
$$\|\Gam(y)\| \leq C (1+\|y\| (1+\|y\|)) G(\tau),$$
where $\lim_{\tau\searrow 0} G(\tau)=0$. 
Then, for any large constant $\bar{K}$, we may choose a sufficiently small 
$\tau$ such that $C (1+\bar{K} (1+\bar{K})) G(\tau)\leq \bar{K}$. 
Hence, we have that, for some small enough 
$\tau<T$,    
$$
\|\Gam(y)\| \leq \bar{K},
\quad\mbox{ whenever }\quad
 \|y\|\leq \bar{K},
$$ 
which obviously means that $\Gam$ maps the closed ball $\mathcal{B}_{\tau,\bar{K}}$ into itself.

\vspace{0.3cm}

We show now that $\Gam$ satisfies a Lipschitz property on $\mathcal{B}_{\tilde{\tau},\bar{K}} \subset \ch^{\ka,\hat{\ka}}(D_{\tilde{\tau}}) $, 
with some $\tilde{\tau}<\tau$ and Lipschitz constant $K<1$. Indeed, owing to the same kind of arguments as before, together with Hypothesis (L), it can be proved that the following estimate holds true:
$$\| \Gam(y_1)-\Gam(y_2)\| \leq C(\tau,\bar{K}) \tilde{G}(\tau) \|y_1-y_2\|,$$
where $\tilde{G}(\tau)$ tends to zero as $\tau$ decreases to zero. Thus, there exists a sufficiently small $\tilde{\tau}<T$ such that $K:=C(\tau,\bar{K}) \|x\|_{\gam,\hat{\gam}} \tilde{G}(\tilde{\tau})<1$. Moreover, we may choose $\tilde{\tau}<\tau$.

\vspace{0.3cm}

Now, all the previous considerations allow us to conclude that $\Gam$ has a unique
fixed point in $\mathcal{H}^{\ka,\hat{\ka}}(D_{\tilde{\tau}})$.

\vspace{0.5cm}

\noindent {\it{Step 2: The {\emph{stairs}} domain $S_{\tilde{\tau}}$ and extension to the square $Q$}}\\
We begin this part of the proof by defining what we understand by the {\it{stairs}} domain $S_{\tilde{\tau}}$. 
First, for a given $(u,v)\in D_T$, we denote by $p_r(u,v)$
the open subset of the line $r:=\{t=-s\} \cap \tilde{R}$ 
corresponding to the projection of the rotated light cone $\tilde{C}(u,v)$
on $r$.

 \vspace{0.3cm}

In step 1 we have proved that Equation (\ref{2.3}) has a unique solution in $\mathcal{H}^{\ka,\hat{\ka}}(D_{\tilde{\tau}})$. In particular, one has existence and uniqueness of solution in the subdomain $S_{\tilde{\tau}} \subset D_{\tilde{\tau}}$, described as follows: 
let $\{P_i=(s_i,t_i), i=1,\dots,N\}$ be a family of points lying on the line $\{t=-s+\sqrt{2}\tilde{\tau}\}$ such that 
$p_r(s_i,t_i)\cap p_r(s_j,t_j)= \emptyset$, if $i\neq j$, and 
$r\cap \tilde{R} \subset \cup_{i=1}^N p_r(s_i,t_i)$.  
Under these conditions, it is clear that $S_{\tilde{\tau}}:= \cup_{i=1}^N \tilde{C}(s_i,t_i)$ forms a 
{\emph{stairs}} domain (see Figure \ref{figure2}). Observe that, since we are considering a finite domain $\tilde{R}$, we may choose a fixed
finite number $N$.

\vspace{0.3cm}

We consider now a square $Q$ determined by two consecutive cones of $S_{\tilde{\tau}}$, say $\tilde{C}(s_i,t_i)$ and $\tilde{C}(s_{i+1},t_{i+1})$, 
for some $i$, as it is shown in Figure \ref{figure2}. In order to simplify notations, we denote by $(a,b)$ the right-upper vertex of $Q$ and $(\hat{a}, b)$, $(a,\hat{b})$ the vertices lying on $\{t=-s+\sqrt{2}\tilde{\tau}\}$, that is $Q=[\hat{a},a]\times[\hat{b},b]$. Set $L_1$ and $L_2$ the sides determined by the points $(\hat{a},\hat{b})$ and $(a,\hat{b})$ and by the former and 
$(\hat{a},b)$, respectively. Notice that we already know the solution on $L_i$, $i=1,2$, and this solution will now play the role
of the initial condition for the equation on the extension $Q$ of  the domain.

\vspace{0.3cm}

We aim to extend the existence and uniqueness result to the square $Q$,
and a first step in this direction, on which we will focus for the remainder
of this step, is to study the invariance of balls in 
$\mathcal{H}^{\ka,\hat{\ka}}(Q)$ under $\Gam$.
The main idea is to decompose the square $Q$ in rectangles of the form $R_n:=[s_n,s_{n+1}]\times [\hat{b},b]$, for $n\geq n_0$, with $n_0$ some positive integer, 
$s_{n_0}=\hat{a}$ and $s_{n+1}-s_n=\frac{1}{n}$ (see Figure \ref{figure2}). Then, we will show that in each space $\mathcal{H}^{\ka,\hat{\ka}}(R_n)$ there is an invariant ball for the corresponding operator $\Gam_n$.  
In the next Step 3 we will focus on the contraction property of the map $\Gam_n$ and thus deduce the existence of a unique fixed point. 
For a given $n\ge n_0$, we will use the notation 
$\|\punt \|_n = \|\punt \|_{\ka,\hat{\ka},n}+ \|\punt \|_{1:\ka,n}+
\|\punt \|_{2:\hat{\ka},n}$ to denote the H\"older semi-norm
$\|\punt \|$ defined for functions on $R_n$.

\vspace{0.3cm}

To begin with, we focus our attention first in the small domain $R_{n_0} \subset Q \subset D_T$. In this case, we are interested in showing the existence of a unique fixed point for the map $\Gam_{n_0}$ defined
for regular functions $y$ and $(s,t)\in R_{n_0}$ as:
$$
(\Gam_{n_0} y)(s,t):= \phi_{n_0}(s,t) + \int \int_{[s_{n_0},s_{n_0+1}]\times [\hat{b},t]} \sig(y(u,v)) x(du,dv).
$$
In the previous relation, the initial condition $\phi_{n_0}$ is given by 
\beq
\phi_{n_0}(s,t)= \psi_{1,n_0}(s,\hat{b}) + \psi_{2,n_0}(s_{n_0},t)+ \psi_{n_0}(s_{n_0},\hat{b}),
\label{q-1}
\eeq
where the first and second term on the right-hand side of the above equality correspond to the known solution on the segments $L_1$ and $L_2$, respectively. The term $\psi_{n_0}(s_{n_0},\hat{b})$ is the initial condition in the corner $(s_{n_0},\hat{b})$ of $R_{n_0}$. 
Notice then that $\|\phi_{n_0}\|<\infty$, which leads us to introduce the 
following set of functions, on which we will build our fixed point argument:
\begin{equation}\label{def:hkno}
\hat\ch^{\ka,\hat\ka}(R_{n_0})=
\lcl
y:R_{n_0}\to\R; \, y_{|_{L_1}}=\psi_{1,n_0}, \, y_{|_{L_2}}=\psi_{2,n_0},
\, \|y\|_{n_0}<\infty
\rcl.
\end{equation}
Observe that, for sake of clarity, we have changed a little the definition
of our functional spaces
with respect to the spaces $\ch^{\ka,\hat\ka}$, by including the initial
condition in the very definition of $\hat\ch^{\ka,\hat\ka}(R_{n_0})$.

\vspace{0.3cm}

Let us fix now a positive number $\mu$ such that $\mu < \ka \land \hat{\ka}$ 
and $\mu< \gam-\ka$. Let also $d$ be a positive constant satisfying 
\beq
\|\phi_{n_0}\|_{n_0} \leq d n_0^\mu.
\label{q0}
\eeq
Our next task is now to prove that the ball $\mathcal{B}_{n_0}$ 
is invariant under $\Gam_{n_0}$, where $\mathcal{B}_{n_0}$ is defined
by:
$$
\mathcal{B}_{n_0}:= \{ y\in \hat{\mathcal{H}}^{\ka,\hat{\ka}}(R_{n_0}); \,  
\|y-\phi_{n_0}\|_{n_0} \leq n_0^\mu \}.$$
We will thus have to
study the norm $\|\Gam_{n_0}y -\phi_{n_0} \|_{n_0}$, and let us 
first consider the H\"older norm $(1:\ka)$.
Let then $(s,t), (s',t)\in R_{n_0}$, with $s'<s$. By the very definition of the operator 
$\Gam_{n_0}$ and Equation (\ref{afitacio}), we have the following estimation
(notice that the initial condition 
$\phi_{n_0}$ cancels out in the first inequality below):
\begin{align}
& |(\Gam_{n_0}y)(s,t)-\phi_{n_0}(s,t) -(\Gam_{n_0}y)(s',t)+\phi_{n_0}(s',t) | \nonumber \\ 
& \quad \leq  
\left| \int \int_{[s',s]\times[\hat{b},t]} \sig(y(u,v)) x(du,dv)\right|\nonumber \\
& \quad \leq   C(s-s')^\gam +C\|y\|_{n_0} (1+\|y\|_{n_0})\left( (s-s')^{\gam+\ka} + (s-s')^\gam \right).
\label{q1}
\end{align}
On the other hand, since $y\in \mathcal{B}_{n_0}$ and we assume condition (\ref{q0}), we have that
\beq
\|y\|_{n_0} \leq \|y-\phi_{n_0}  \|_{n_0} + \|\phi_{n_0}  \|_{n_0} \leq (d+1) n_0^\mu.
\label{q1.5}
\eeq
Plugging this bound into (\ref{q1}) and taking into account that $s_{n_0+1}-s_{n_0}=\frac{1}{n_0}$, it turns out that
\beq
\|\Gam_{n_0}y -\phi_{n_0}\|_{1:\ka,n_0} \leq  \frac{1}{n_0^{\gam-\ka}}\big(1+ (d+1)n_0^\mu (1+ (d+1)n_0^\mu)\big)
=o(n_0^\mu).
\label{q2}
\eeq
Analogously, we obtain the following estimate for the H\"older norm  $\| \Gam_{n_0}y -\phi_{n_0}\|_{2:\hat{\ka},n_0}$:
\beq
\|\Gam_{n_0}y -\phi_{n_0}\|_{2:\hat{\ka},n_0} \leq  \frac{1}{n_0^\gam} \big(1+ (d+1) n_0^\mu(1+ (d+1)n_0^\mu)\big)
=o(n_0^\mu).
\label{q3}
\eeq
Eventually, let us deal with the H\"older norm on rectangles. Let $R=[s_1,s_2]\times [t_1,t_2]$ be a rectangle included in $R_{n_0}$. It is readly checked that $\Del_R \phi_{n_0}=0$. Thus, owing to (\ref{afitacio}), 
we have that
\begin{align*}
|\Del_R (\Gam_{n_0}y) |\leq & C(s_2-s_1)^\gam (t_2-t_1)^{\hat{\gam}} \\
&\quad + C \|y\|_{n_0} (1+\|y\|_{n_0}) \big( (s_2-s_1)^{\gam+\ka} (t_2-t_1)^{\hat{\gam}}+ (s_2-s_1)^\gam (t_2-t_1)^{\hat{\gam}+\hat{\ka}}\big).
\end{align*}
Hence, by (\ref{q0}), (\ref{q1.5})  and the fact that $y\in \mathcal{B}_{n_0}$, we end up with
\beq
\|\Gam_{n_0}y -\phi_{n_0}\|_{\ka,\hat{\ka},n_0} \leq  \frac{C}{n_0^{\gam -\ka}} \big( 1+ (d+1)n_0^\mu(1+ (d+1)n_0^\mu)\big)= o(n_0^\mu). 
\label{q4}
\eeq
Putting together (\ref{q2})-(\ref{q4}), we obtain that the ball $\mathcal{B}_{n_0}$ is invariant under $\Gam_{n_0}$ for $n_0$ large enough.

\vspace{0.3cm}

At this point, let us anticipate a little on the next step, and assume
our contraction arguments have lead us to the definition of a unique
solution up to the rectangle $R_n$ such that $R_n\cap Q\neq \emptyset$,
for $n\ge n_0$.
We will then try to use an induction argument in order to define
an invariant ball under the map $\Gam_{n+1}$.
Since $s_{n+1}-s_n=\frac{1}{n}$, even if the size of $R_n$ 
decreases, we will cover the whole square $Q$ in a finite number of
steps, thanks to the fact that $\sum n^{-1}$ is a divergent series.
Observe then that, if the solution $y$ to our equation has been defined up to
$R_n$, and if $L_2^{n+1}$ denotes the left vertical side of $R_{n+1}$,
then the solution to  (\ref{2.3}) on $R_{n+1}$ should satisfy
$$
y_{|L_1}=\psi_{1,n+1}, 
\quad\mbox{ and }\quad
y_{|L^{n+1}_2}=\psi_{2,n+1},
$$
for the function $\psi_{2,n+1}=y_{|R_n\cap L^{n+1}_2}$, and where 
$\psi_{1,n+1}$ has been introduced at relation (\ref{q-1}).
We will thus introduce a space $\hat\ch^{\ka,\hat\ka}(R_{n+1})$ analogously
to the case $n=n_0$ given at (\ref{def:hkno}):
$$
\hat\ch^{\ka,\hat\ka}(R_{n+1})=
\lcl
y:R_{n+1}\to\R; \, y_{|L_1\cap R_{n+1}}=\psi_{1,n+1},\, 
y_{|L_2^{n+1}}=\psi_{2,n+1},\,
\|y\|_{n+1}<\infty
\rcl.
$$
Assume now that, for any $n_0\le m\le n$, the operator 
$\Gam_m:\hat\ch^{\ka,\hat\ka}(R_m)\to\hat\ch^{\ka,\hat\ka}(R_m)$, 
defined by
$$
(\Gam_{m} y)(s,t):= \phi_m(s,t) + \int \int_{[s_m,s_{m+1}]
\times [\hat{b},t]} \sig(y(u,v)) x(du,dv),$$
for $y\in \mathcal{H}^{\ka,\hat{\ka}}(R_m)$ and $(s,t)\in R_m$, leave the 
following ball invariant:
$$
\mathcal{B}_m= \{ y\in \hat\ch^{\ka,\hat{\ka}}(R_m); \|y-\phi_n\|_n \leq n^\mu \},
$$
where $\phi_m$ is defined as in 
(\ref{q-1}). We aim to show that the same is true on $R_{n+1}$. 
This can be achieved using the same kind of calculations 
as for the case $n=n_0$ and applying the following result:
\begin{lemma}
For all $n\geq n_0$, it holds that $$\|\phi_n\|_n \leq d n^\mu.$$
\label{lemmaq}
\end{lemma}
\noindent {\it{Proof:}} Our statement holds true for $n=n_0$, by hypothesis. 
Assume that we have proved, for some $n>n_0$, that 
$\|\phi_k\|_k \leq d k^\mu$ for any $k\le n$, and let us prove the
property for $\phi_{n+1}$. 

\vspace{0.3cm}

Firstly, it is straightforward to check that, for all $k\geq n_0$, the function $\phi_k$ has null rectangular increments, namely 
$\|\phi_k\|_{\ka,\hat{\ka},k}=0$. Moreover, it holds that $\|\phi_{n+1}\|_{1:\ka,n+1}\leq  \|\psi_{1,n+1}\|_{1:\ka,n+1}$. 
On the other hand,
observe that the side $L_1$ does not vary when we jump from $n$ to $n+1$, for all $n\geq n_0$, 
which implies that $\psi_{1,n+1}=\psi_{1,n_0}$. Hence, the norm 
$\|\psi_{1,n+1}\|_{1:\ka,n+1}$ may be bounded by a constant independent of $n$, say $C_0$.

\vspace{0.3cm}

Let us seek now some estimates for the quantity $\|\phi_{2,n+1}\|_{2:\hat{\ka},n+1}$. Let $s\in (s_{n+1}, s_{n+2})$ and 
$t,t' \in (\hat{b},b)$. Then, by the very definition of $\phi_{n+1}$, we have the following estimation:
\begin{align*}
|\phi_{n+1}(s,t)-\phi_{n+1}(s,t')| & =
|\psi_{2,n+1}(s_{n+1},t)-\psi_{2,n+1}(s_{n+1},t')| \\
& \leq  |\Del_{U_n} \psi_{2,n+1}| 
+ | \psi_{2,n}(s_n,t)-\psi_{2,n}(s_n,t')|\\
& \leq \Big( \| y \|_{\ka,\hat{\ka},n} (s_{n+1}-s_n)^\ka  + \|\psi_{2,n}\|_{2:\hat{\ka},n} \Big) (t-t')^{\hat{\ka}},
\end{align*}
where $U_n:=[s_n,s_{n+1}]\times[t',t]$ and $y$ denotes the unique 
solution to our equation in $\hat{\mathcal{H}}^{\ka,\hat{\ka}}(R_n)$. Therefore, making use of the induction hypothesis, we end up with
$$\|\psi_{2,n+1}\|_{2:\hat{\ka},n+1} \leq \frac{1}{n^{\ka-\mu}} + d n^\mu .$$
Thus, we have proved that 
$$\|\phi_{n+1}\|_{n+1} \leq C_0+ \frac{1}{n^{\ka-\mu}} + d n^\mu \leq d(n+1)^\mu,$$
which concludes the proof of the lemma. \qed

\vspace{0.3cm}

Summing up, we have defined a finite sequence of rectangles $(R_n)_n$ such that $Q\subset \cup_{n} R_n$ and we have proved that there exists a ball 
$\cb_n$ in 
$\hat\ch^{\ka,\hat{\ka}}(R_n)$ which is invariant under $\Gam_n$. In the next Step 3 we shall address the contractivity properties of $\Gam_n$.
Notice again that we have assumed in advance that the equation could
be solved on any of the rectangles $R_k$ for $k\le n$ in order to
define the ball $\cb_{n+1}$. This claim will be justified at Step 3.

\vspace{0.3cm}

\noindent {\it{Step 3: Contraction property}}\\
In this part of the proof, we fix $R_n$ any of rectangles covering the square $Q$ and we consider $\mathcal{B}_n$ the invariant ball for $\Gam_n$ in 
$\hat\ch^{\ka,\hat{\ka}}(R_n)$. We aim to prove that $\Gam_n$ has a unique fixed point in $\mathcal{B}_n$. 

\vspace{0.3cm}

In order to show that $\Gam_n$ is a contraction, we will consider again a suitable family of squares $(\theta_j, j=1,\dots,M)$ covering $R_n$ and 
having the same size. Then we will prove that the operator $\Gam_n$, restricted to any $\theta_j$, is a contraction, that is, $\Gam_n$ defined on $\mathcal{H}^{\ka,\hat{\ka}}(\theta_j)$ satisfies the following condition:
\beq
\| \Gam_n(y_1)-\Gam_n(y_2)\|_{n, \theta_j} \leq K \| y_1-y_2\|_{n, \theta_j},
\label{q5}
\eeq
for all $y_1, y_2 \in \mathcal{H}^{\ka,\hat{\ka}}(\theta_j)$ and for some constant $K<1$;
$\| \punt\|_{n, \theta_j}$ denotes the corresponding norm on $\mathcal{H}^{\ka,\hat{\ka}}(\theta_j)$.  

\vspace{0.3cm}

Recall that $R_n=[s_n,s_{n+1}]\times [\hat{b},b]$. Let $\theta_1$ be the square contained in $R_n$ defined by 
$[s_n, \tilde{a}]\times [\hat{b}, \tilde{b}]$, for some $\tilde{a}, \tilde{b}$ (See Figure \ref{figure2}). In this case, we are interested in the 
operator $\Gam_{n,\theta_1}$, defined on $\mathcal{H}^{\ka,\hat{\ka}}(\theta_1)$, as follows:
$$(\Gam_{n,\theta_1} y)(s,t)= \phi_{n,\theta_1}(s,t) + \int \int_{[s_n,s]\times[\hat{b},t]} \sig(y(u,v)) x(du,dv),$$
$(s,t)\in \theta_1$,  where $\phi_{n,\theta_1}$ corresponds to the initial condition. Let us prove that if $\theta_1$ has a sufficiently small size, then 
$\Gam_{n,\theta_1}$ is a contraction.

\vspace{0.3cm}

To begin with, let us deal with the H\"older norm in rectangles. Namely, let $y_1, y_2 \in \mathcal{H}^{\ka,\hat{\ka}}(\theta_1)$ and $R$ a 
rectangle contained in $\theta_1$, so we study the following expression:
$$\left| \Del_R \left(\Gam_{n,\theta_1} (y_1) - \Gam_{n,\theta_1} (y_2)\right) \right| = 
\left| \int \int_R (\sig(y_1(u,v))- \sig(y_2(u,v))) x(du,dv)\right|.$$
Owing to (\ref{afitacio}), Hypothesis (L), the fact that $y_1, y_2\in \mathcal{B}_n$ and Lemma \ref{lemmaq}, one has the following estimate:
\beq
\| \Gam_{n,\theta_1} (y_1) - \Gam_{n,\theta_1} (y_2)\|_{\ka,\hat{\ka},\theta_1} \leq \|y_1-y_2\|_{n,\theta_1} C(n,\mu) 
(\tilde{a}-s_n)^{\gam-\ka} (\tilde{b}-\hat{b})^{\hat{\gam}-\hat{\ka}}.
\label{q6}
\eeq
Concerning the H\"older norms $(1:\ka)$ and $(2:\hat{\ka})$, one uses similar arguments as for the above norm on rectangles to end up with bounds for 
$\| \Gam_{n,\theta_1} (y_1) - \Gam_{n,\theta_1} (y_2)\|_{1:\ka,\theta_1}$ and  
$\| \Gam_{n,\theta_1} (y_1) - \Gam_{n,\theta_1} (y_2)\|_{2:\hat{\ka},\theta_1}$ of the same type as (\ref{q6}). Therefore,   
if the size of $\theta_1$, say $\del:=\tilde{a}-s_n$, is sufficiently small,  
then we have that condition (\ref{q5}), for $j=1$, is fulfilled.  

\vspace{0.3cm}

Arguing as for the square $\theta_1$, one could see that we can progressively cover the rectangle $R_n$ 
by a finite family of squares $(\theta_j, j=1,\dots,M)$, such that each $\theta_j$ has the same size $\del$ and the corresponding operator 
$\Gam_{n,\theta_j}$ is a contraction on  
$\mathcal{H}^{\ka,\hat{\ka}}(\theta_j)$ (notice however that $\delta$ depends 
on $n$). 
Thus, this let us conclude that the 
map $\Gam_n : \mathcal{H}^{\ka,\hat{\ka}}(R_n)\rightarrow \mathcal{H}^{\ka,\hat{\ka}}(R_n)$ has a unique fixed point. 

\vspace{0.5cm}

 \noindent {\it{Step 4: Extension to the whole domain}}\\
Putting together the considerations of Step 2 and Step 3, we have 
constructed now, with a finite number of steps, a unique solution to
(\ref{2.3}) on the whole square $Q$.
Analogously, we will be able to obtain the same result for 
the other squares determined by the stairs domain $S_{\tilde{\tau}}$, in such a way that we have proved the existence and uniqueness of solution to Equation (\ref{2.3}) in an extended {\it{stairs}} domain $S_{\tau_1}$, 
for some $\tau_1 > \tilde{\tau}$. Eventually, we iterate this procedure in order to cover the whole 
domain $D_T$ of definition of our Equation (\ref{2.3}). This concludes the 
proof of the theorem.

\qed

\vspace{0.3cm}

\subsection{Application to the fractional Brownian motion}\label{applifbm}

In this section we apply Theorem \ref{theorem1} in the particular case where the function $x$ corresponds to the path of some random perturbation. Namely, we are interested in Gaussian random noises having a fractional time correlation and some spatially homogeneous one. Let us make this rigorous, as follows.

\vspace{0.3cm}

Fix $H\in (\frac{1}{2},1)$ and consider,
on a given complete probability space $(\Omega,\cf,P)$,
 a $L^2(\Om) -$valued centered Gaussian process $\{X(\phi), \phi\in \mathcal{D}(\re^2)\}$, where $\mathcal{D}(\re^2)$ stands for the space of test functions, with covariance functional given by 
$$E\left(X(\phi)X(\psi)\right)= c_H \int_{[0,T]^2} du dv |u-v|^{2H-2} \int_{\re^2}dx dy Q(x-y) \phi(u,x)\psi(v,y).$$
The value of $c_H$ is $H(2H-1)$ and we will focus on the case in which $Q$ is a Riesz kernel, that is 
$Q(x)=\frac{1}{|x|^\nu}$, for some $\nu \in (0,1)$.

\vspace{0.3cm}

This section will be decomposed as follows:
first, we will define a process denoted by $\{x(s,t), (s,t)\in D_T\}$, which will correspond to the 45º rotation of $X$ and will be the driving motion of Equation (\ref{2.3}). Then, we will study the regularity properties of $x$, namely the Hölder regularity in rectangles, and therefore deduce that Theorem \ref{theorem1} can be applied, 
which will give a
proof of Theorem \ref{fbmsol}.
Notice that,
in order to obtain the $\mathcal{C}^{\gam,\hat{\gam}}-$regularity of $x$, for some $\gam,\hat{\gam}$, we will apply the following extension of Kolmogorov's lemma (see \cite{feyel}).
\begin{lemma}
Let $z=\{z(s,t), (s,t)\in \re^2\}$ be a two-parameter stochastic process. Assume that there exist $p>1$ and  $a,b\in (\frac{1}{p},+\infty)$ such that
$$\|\Del_R z\|_{L^p(\Om)} \leq C (s_2-s_1)^a (t_2-t_1)^b,$$
for any rectangle $R=[s_1,s_2]\times [t_1,t_2]$ and some positive constant $C$. Then, the process $z$ admits a continuous modification whose trajectories belong to the space $\mathcal{C}^{a'-\frac{1}{p},b'-\frac{1}{p}}$, 
for all $a'\in (\frac{1}{p},a)$ and $b'\in (\frac{1}{p},b)$. 
\label{lemma2}
\end{lemma}

\vspace{0.3cm}

Let us begin now with the construction of the rotation of $X$:
we denote by $\mathcal{R}_{\frac{\pi}{4}}$ the 45º rotation on the plane,
and for any function $\ffi\in \mathcal{D}(D_T)$, we set 
$$\tilde{X}(\ffi):=X(\ffi \circ  \mathcal{R}_{\frac{\pi}{4}}).$$
It is straighforward to check that $\{ \tilde{X}(\ffi), \ffi\in \mathcal{D}(D_T)\}$ is a well defined $L^2(\Om)-$valued centered Gaussian process with covariance functional
\begin{align*}
&E\left(\tilde{X}(\ffi_1)\tilde{X}(\ffi_2)\right)=  c_H \int_{[0,T]^2} du dv |u-v|^{2H-2}\\
& \; \times  \int_{\re^2}dx dy |x-y|^{-\nu} \ffi_1\left(\frac{x-u}{\sqrt{2}},\frac{x+u}{\sqrt{2}}\right)
\ffi_2\left(\frac{y-v}{\sqrt{2}},\frac{y+v}{\sqrt{2}}\right).
\end{align*}
Now, we aim to define $\tilde{X}(\1_R)$, where $R$ is some rectangle contained in $D_T$. For this, we 
consider a sequence of functions $(\ffi_n)_n$ in $\mathcal{D}(D_T)$ such that $\ffi_n(s,t)\rightarrow \1_R (s,t)$ as $n$ tends to infinity, for all $(s,t)\in D_T$. Then, by bounded convergence, it follows that
$$\lim_{n,m\rightarrow \infty} E(|\tilde{X}(\ffi_n)-\tilde{X}(\ffi_m)|^2)=   
\lim_{n,m\rightarrow \infty} E(|\tilde{X}(\ffi_n-\ffi_m)|^2)=0,$$
and we set 
$$\tilde{X}(\1_R):=L^2(\Om)-\lim_{n\rightarrow \infty} \tilde{X}(\ffi_n).$$
It can be easily seen that this definition does not depend on the particular approximating sequence. Moreover, we still have the following equality,
for two arbitrary rectangles $R_1,R_2$ of $D_T$:
\begin{align}
&E\left(\tilde{X}(\1_{R_1})\tilde{X}(\1_{R_2})\right)=  c_H \int_{[0,T]^2} du dv |u-v|^{2H-2}\nonumber \\
& \; \times  \int_{\re^2}dx dy |x-y|^{-\nu} \1_{R_1}\left(\frac{x-u}{\sqrt{2}},\frac{x+u}{\sqrt{2}}\right)
\1_{R_2}\left(\frac{y-v}{\sqrt{2}},\frac{y+v}{\sqrt{2}}\right).
\label{4.1}
\end{align}
It is then natural to define, for $(s,t)\in D_T$, 
\begin{equation}\label{defxfbm}
x(s,t):=\tilde{X}\left( \1_{\left([0,s]\times[0,t]\right)\cap D_T}\right),
\end{equation}
with the convention that, if either $s$ or $t$ is negative, the corresponding interval will be $[s,0]$ or $[t,0]$, 
respectively. Recall that $x$ will be the driving noise of equation
(\ref{2.3}).

Let us fix an element $\om$ of the probability space $(\Om, \tf,P)$ on which the process $X$ is defined. In order to simplify the notation, we shall still denote by $\{x(s,t), (s,t)\in D_T\}$ the path of the process $x$ associated to $\om$. Then, the following result caracterises the Hölder regularity on rectangles of the function $x$, which is obviously an important step in 
order to apply Theorem \ref{theorem1}.
\begin{lemma}
Almost surely on $\Om$,
the function  $x:D_T\rightarrow \re$ defined by (\ref{defxfbm})
belongs to the space $\mathcal{C}^{\eta,\hat{\eta}}$, with 
$\eta\in (0,\gam)$ and $\hat{\eta}\in (0,\hat{\gam})$, for any $\gam,\hat{\gam}$ such that $\gam+\hat{\gam}=H+ \frac{2-\nu}{2}$. 
\label{lemma3}
\end{lemma}  

\vspace{0.5cm}

\noindent
{\it{Proof}}. We will apply Lemma \ref{lemma2}. For this, we fix a rectangle $R=[s_1,s_2]\times[t_1,t_2]$ in $D_T$ and we compute the square moment of $\Del_R x$. Indeed, by the very definition of $x$, we obtain that
$$E(|\Del_R x|^2)=E(|\tilde{X}(\1_R)|^2)=E(|X(\1_R \circ R_{\frac{\pi}{4}})|^2).$$
Hence, owing to (\ref{4.1}), we end up with 
\begin{align}
E\left(|\Del_R x|^2\right)= & c_H \int_{[0,T]^2} du dv |u-v|^{2H-2}\nonumber \\
& \; \times  \int_{\re^2}dx dy |x-y|^{-\nu} \1_R \left(\frac{x-u}{\sqrt{2}},\frac{x+u}{\sqrt{2}}\right)
\1_R \left(\frac{y-v}{\sqrt{2}},\frac{y+v}{\sqrt{2}}\right),
\label{4.2}
\end{align}  
and notice that, for instance, we have 
$\1_R \left(\frac{x-u}{\sqrt{2}},\frac{x+u}{\sqrt{2}}\right)=
\1_{\hat{R}}(u,x)$, where $\hat{R}:=\mathcal{R}_{-\frac{\pi}{4}}(R)$.

\vspace{0.3cm}

In the sequel, we will make use of the following notation:
$$S_{ij}:=\frac{t_i+s_j}{\sqrt{2}},\;\; T_{ij}:=\frac{t_i-s_j}{\sqrt{2}},\;\; i,j=1,2, 
$$
and set also $\Del s:=s_2-s_1$ and $\Del t:=t_2-t_1$.
With these notations in hand, we will try to get a bound of the type
\begin{equation}\label{incx}
E\left( |\Del_R x|^2\right) \leq C (\Del s)^{2\gam}(\Del t)^{2\hat{\gam}},
\end{equation}
where $\gam,\hat\gam$ satisfy $\gam+\hat{\gam}=H+ \frac{2-\nu}{2}$,
for a general rectangle $R$.
However, throughout the proof, we will assume that 
$\Del s < \Del t$: indeed the case  $\Del t < \Del s$ may be treated 
analogously, 
and when $R$ is a square, relation (\ref{incx}) follows easily.

\vspace{0.3cm}

\begin{figure}[!ht]
\begin{center}
\scalebox{0.40}{\includegraphics{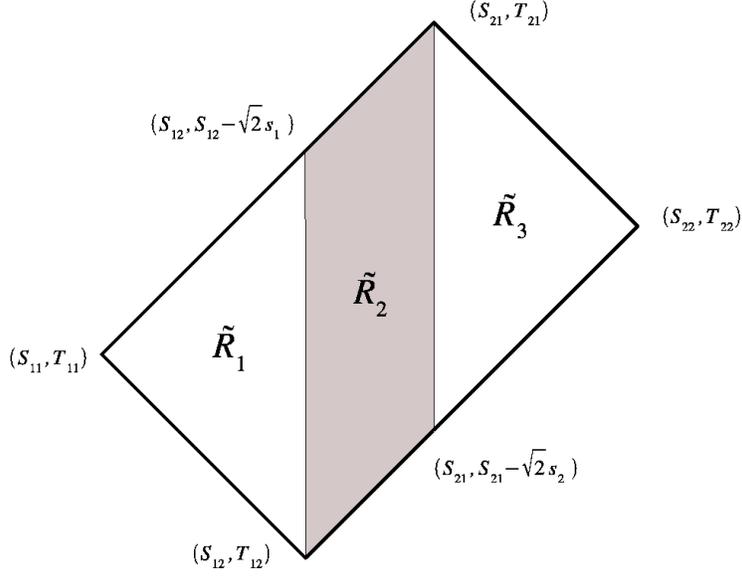}}   
\end{center}
\caption{Decomposition of the domain $\hat{R}=\mathcal{R}_{-\frac{\pi}{4}}(R)$.}
\label{figure3}
\end{figure}

In order to get good bounds of the right-hand side of (\ref{4.2}), we will decompose the indicator function
$\1_{\hat{R}}$, as follows (see Figure \ref{figure3}):
$$
\1_{\hat{R}}=\1_{\tilde{R}_1} +\1_{\tilde{R}_2} +\1_{\tilde{R_3}},
$$
where 
\begin{itemize}
\item $\tilde{R}_1$ is the triangle determined by
$$\{(S_{11},T_{11}),(S_{12},T_{12}),(S_{12},S_{12}-\sqrt{2}s_1)\},$$
\item $\tilde{R}_2$ is the parallelogram determined by
$$\{(S_{12},T_{12}),(S_{12},S_{12}-\sqrt{2}s_1),(S_{21},T_{21}),(S_{21},S_{21}-\sqrt{2}s_2)\},$$
\item $\tilde{R}_3$ is the triangle determined by
$$\{(S_{21},T_{21}),(S_{21},S_{21}-\sqrt{2}s_2),(S_{22},T_{22})\}.$$
\end{itemize} 
Owing to this decomposition, the integral in the right-hand side of (\ref{4.2}) shall be split in six terms.

\vspace{0.3cm}

The proof can be divided in two different situations:
 
\noindent {\it{First case: $\Del t < 3\Del s.$}} Under this assumption, it is straightforward to check that, in the terms of the integral (\ref{4.2}) involving the indicator function of the parallelogram  $\tilde{R}_2$ -which are the most winding-, there will always be intersection between the domains of variation of $x$ and $y$. Moreover, the fact that $\Del s< \Del t < 3\Del s$ will let us obtain the appropiate bounds for (\ref{4.2}) without too much effort. The details for this case are omitted here.

\noindent {\it{Second case: $\Del t \geq 3\Del s.$}} Here, fixed $u$ and $v$ in their respective domains of variation, there does not need to be intersection between the domains of variation of $x$ and $y$ in the terms involving the parallelogram $\tilde{R}_2$. As it has been mentioned, it turns out that this situation leads to much more complicated calculations in comparison with the above one. Therefore, we will just deal with the proof in this second case.

\vspace{0.3cm}

As mentioned before, our partition of $\hat{R}$ into three pieces leads
to the computation of six different terms. However,
by symmetry, it is sufficient to consider the following three situations:   

\vspace{0.3cm}

\begin{center}
{\underline{Term corresponding to $\1_{\tilde{R}_2}(u,x) \1_{\tilde{R}_2}(v,y)$}}
\end{center}
We have to get bounds, for instance, of the integral 
\beq
\int_{S_{12}}^{S_{21}} du\int_{S_{12}}^u dv \int_{u-\sqrt{2}s_2}^{u-\sqrt{2}s_1} dx 
\int_{v-\sqrt{2}s_2}^{v-\sqrt{2}s_1} dy (u-v)^{2H-2} |x-y|^{-\nu};
\label{7}
\eeq
the case when $v\in (u,S_{21})$ can be treated in the same way. The main idea is to decompose the domain of integration above so as to get rid of the absolute value there, and therefore be able to compute, or eventually bound, the remaining integrals. Notice that, for the sequel of this proof, we have decided not to write the multiplicative constants which appear in our estimations, and which can vary from line to line.

\vspace{0.3cm}  

First, it holds that (\ref{7}) equals to
\begin{align*}
&A_1+A_2 :=   
\int_{S_{12}}^{S_{12}+\sqrt{2}\Del s} du\int_{S_{12}}^u dv \int_{u-\sqrt{2}s_2}^{u-\sqrt{2}s_1} dx 
\int_{v-\sqrt{2}s_2}^{v-\sqrt{2}s_1} dy (u-v)^{2H-2} |x-y|^{-\nu}\\
&+ \int_{S_{12}+\sqrt{2}\Del s}^{S_{21}} du\int_{S_{12}}^u dv \int_{u-\sqrt{2}s_2}^{u-\sqrt{2}s_1} dx 
\int_{v-\sqrt{2}s_2}^{v-\sqrt{2}s_1} dy (u-v)^{2H-2} |x-y|^{-\nu}.
\end{align*}
Then, on one hand we have that
\begin{align*}
A_1= & \int_{S_{12}}^{S_{12}+\sqrt{2}\Del s} du\int_{S_{12}}^u dv \int_{u-\sqrt{2}s_2}^{v-\sqrt{2}s_1} dx 
\int_{v-\sqrt{2}s_2}^x dy (u-v)^{2H-2} (x-y)^{-\nu}\\
&+ \int_{S_{12}}^{S_{12}+\sqrt{2}\Del s} du\int_{S_{12}}^u dv \int_{u-\sqrt{2}s_2}^{v-\sqrt{2}s_1} dx 
\int_x^{v-\sqrt{2}s_1} dy (u-v)^{2H-2} (y-x)^{-\nu}\\
&+ \int_{S_{12}}^{S_{12}+\sqrt{2}\Del s} du\int_{S_{12}}^u dv \int_{v-\sqrt{2}s_1}^{u-\sqrt{2}s_1} dx 
\int_{v-\sqrt{2}s_2}^{v-\sqrt{2}s_1} dy (u-v)^{2H-2} (x-y)^{-\nu},
\end{align*}
and it is readily checked,
from this decomposition, that 
$$A_1\leq  (\Del s)^{2H-\nu+2}.$$
On the other hand, one has that $A_2=B_1+B_2$, with
\begin{align*}
B_1:= &  
\int_{S_{12}+\sqrt{2}\Del s}^{S_{21}} du \int_{S_{12}}^{u-\sqrt{2}\Del s} dv \int_{u-\sqrt{2}s_2}^{u-\sqrt{2}s_1} dx 
\int_{v-\sqrt{2}s_2}^{v-\sqrt{2}s_1} dy (u-v)^{2H-2} (x-y)^{-\nu},\\
B_2:= & 
\int_{S_{12}+\sqrt{2}\Del s}^{S_{21}} du \int_{u-\sqrt{2}\Del s}^u dv \int_{u-\sqrt{2}s_2}^{u-\sqrt{2}s_1} dx 
\int_{v-\sqrt{2}s_2}^{v-\sqrt{2}s_1} dy (u-v)^{2H-2} |x-y|^{-\nu}.
\end{align*}
Notice that, in the domain of integration of $B_1$,
the following estimation holds true:
$$
(x-y)^{-\nu}\leq (u-v-\sqrt{2}\Del s)^{-\nu}.
$$
Hence, 
$$B_1 \leq  (\Del s)^2  
\int_{S_{12}+\sqrt{2}\Del s}^{S_{21}} du \int_{S_{12}}^{u-\sqrt{2}\Del s} dv  (u-v)^{2H-2},$$
and computing the above integral, one ends up with the bound
$$B_1 \leq  (\Del s)^2 (\Del s + \Del t)^{2H-\nu}.$$
Eventually, it turns out that we can decompose $B_2$ as follows:
\begin{align*}
B_2= & \int_{S_{12}+\sqrt{2}\Del s}^{S_{21}} du \int_{u-\sqrt{2}\Del s}^u dv
 \int_{u-\sqrt{2}s_2}^{v-\sqrt{2}s_1} dx 
\int_{v-\sqrt{2}s_2}^x dy (u-v)^{2H-2} (x-y)^{-\nu} \\
&+  \int_{S_{12}+\sqrt{2}\Del s}^{S_{21}} du \int_{u-\sqrt{2}\Del s}^u dv\int_{u-\sqrt{2}s_2}^{v-\sqrt{2}s_1} dx 
\int_x^{v-\sqrt{2}s_1} dy (u-v)^{2H-2} (y-x)^{-\nu} \\
&+  \int_{S_{12}+\sqrt{2}\Del s}^{S_{21}} du \int_{u-\sqrt{2}\Del s}^u dv\int_{v-\sqrt{2}s_1}^{u-\sqrt{2}s_1} dx 
\int_{v-\sqrt{2}s_2}^{v-\sqrt{2}s_1} dy (u-v)^{2H-2} (x-y)^{-\nu},
\end{align*}
and from this expression, it is not difficult to prove that
$$B_2\leq (\Del s)^{2H-\nu+1} (\Del t + \Del s).$$
Putting together all the bounds that we have obtained on
$B_1$, $B_2$ $A_2$ and $A_1$, we conclude that the term (\ref{7}) 
can be bounded by
\beq
(\Del s)^{2H-\nu+2} + (\Del s)^2 (\Del s + \Del t)^{2H-\nu} + (\Del s)^{2H-\nu+1} (\Del t + \Del s).
\label{4.4}
\eeq
Let $\gam$, $\hat{\gam}$ belong to $(\frac{1}{2},1)$ and such that $\gam +\hat{\gam}=H+\frac{2-\nu}{2}$. Then, owing to (\ref{4.4}), one can easily see that (\ref{7}) may be bounded, up to constants, by $(\Del s)^{2\gam}(\Del t)^{2\hat{\gam}}$.

\vspace{0.3cm}

\begin{center}
{\underline{Term corresponding to $\1_{\tilde{R}_1}(u,x) \1_{\tilde{R}_2}(v,y)$}}
\end{center}
Let us treat now
the term corresponding to $\1_{\tilde{R}_1}(u,x) \1_{\tilde{R}_2}(v,y)$, namely\beq
\int_{S_{11}}^{S_{12}} du\int_{S_{12}}^{S_{21}} dv \int_{-u+\sqrt{2}t_1}^{u-\sqrt{2}s_1} dx 
\int_{v-\sqrt{2}s_2}^{v-\sqrt{2}s_1} dy (v-u)^{2H-2} |x-y|^{-\nu}.
\label{4.5b}
\eeq
This term equals to $D_1 +D_2$, where
\begin{align*}
D_1 = & \int_{S_{11}}^{S_{12}} du\int_{S_{12}}^{S_{12}+\sqrt{2}\Del s} dv 
 \int_{-u+\sqrt{2}t_1}^{u-\sqrt{2}s_1} dx 
\int_{v-\sqrt{2}s_2}^{v-\sqrt{2}s_1} dy
(v-u)^{2H-2} |x-y|^{-\nu},\\
D_2 = & \int_{S_{11}}^{S_{12}} du\int_{S_{12}+\sqrt{2}\Del s}^{S_{21}} dv 
 \int_{-u+\sqrt{2}t_1}^{u-\sqrt{2}s_1} dx 
\int_{v-\sqrt{2}s_2}^{v-\sqrt{2}s_1} dy
(v-u)^{2H-2} (y-x)^{-\nu},
\end{align*}
and using similar arguments as for the term $A_1$, one can easily 
show that the following estimate holds true 
(recall that we omitt to write the 
multiplicative constants):
$$D_2 \leq (\Del s)^2 (\Del t)^{2H-\nu}.$$
For the term $D_1$, we can show, along the 
 same lines as for the term $A_2$ above, that
$$
D_1 \leq (\Del s)^{2H-\nu+2}.
$$
Thus, (\ref{4.5b}) may be estimated by 
$$(\Del s)^{2\gam}(\Del t)^{2\hat{\gam}},$$
with $\gam +\hat{\gam}=H+\frac{2-\nu}{2}$. 

\vspace{0.3cm}

\begin{center}
{\underline{Term corresponding to $\1_{\tilde{R}_1}(u,x) \1_{\tilde{R}_1}(v,y)$}}
\end{center}
Our aim now is to bound the integral
$$\int_{S_{11}}^{S_{12}} du\int_{S_{11}}^{S_{12}} dv \int_{-u+\sqrt{2}t_1}^{u-\sqrt{2}s_1} dx 
\int_{-v+\sqrt{2}t_1}^{v-\sqrt{2}s_1} dy |u-v|^{2H-2} |x-y|^{-\nu}.
$$
It turns out that this term involves easier calculations in comparison with the two other situations, and that 
it can be bounded by 
$$(\Del s)^{2H-\nu +2}\leq (\Del s)^{2\gam}(\Del t)^{2\hat{\gam}}.$$
The details of computations for this term are left to the reader.

\vspace{0.3cm}

Putting together the estimations obtained in the study of the three terms
 above, we get that
$$E\left( |\Del_R x|^2\right) \leq C (\Del s)^{2\gam}(\Del t)^{2\hat{\gam}},$$
for a positive constant $C$.
Taking into account that we are in a Gaussian context and applying Lemma \ref{lemma2}, we get the desired regularity for the process $x$. \qed

\begin{remark}
As mentioned before, Theorem \ref{fbmsol} is now an easy consequence
of the previous lemma. Indeed, it allows us to 
apply Theorem \ref{theorem1} to the fractional Brownian motion $x$
defined by (\ref{defxfbm}), 
which in turn implies  Theorem \ref{fbmsol},
by the simple rotation trick of Section \ref{rotwa}.
\end{remark}

\section{Sharpness of our method}\label{optim}

Let us come back for a moment to Section \ref{rotwa}, and to the definition
(\ref{2.3}) we gave for the wave equation: we performed a rotation
$\mathcal{R}_{\frac{\pi}{4}}$ on our initial equation in order to get a new domain
of integration $\tilde{C}(s,t)$ whose sides are parallel to the axes.
Then, as pointed out in Remark \ref{remark2}, if $x,y:D_T \rightarrow \re$ 
are two functions satisfying Hypothesis (H) and $\sig:\re \rightarrow \re$ 
fulfills Hypothesis (L), the integral 
$$I(s,t)\equiv\int \int_{\tilde{C}(s,t)} \sig(y(u,v)) x(du,dv)$$ 
is well defined in the Young 
sense. Moreover, it has the same H\"older regularity on rectangles as the 
function $x$. 

However, as it will be made clearer later on, one could define directly the integral
$$
\mathcal{I}(s,t)=\int \int_{C(s,t)} \sig(Y(u,v))X(du,dv)
=\int_0^s \int_{\re} G_{s-u}(t,v) \sig(Y(u,v))X(du,dv),
$$
where the fundamental solution $G$ of the wave equation is 
given by
$$G_s(t,v)=\frac{1}{2} \1_{ \{|t-v|<s\} },\; s>0, \, t,v\in \re,$$
by means of Young approximations with respect to a {\emph{good}} partition of the domain, such as a dyadic one. 
In this section, we will try to show that this direct strategy does not 
behave as well as the one we proposed at Section \ref{weq}. Indeed,
we will show that, even in the linear case (i.e. $\sig\equiv 1$), we are quite far away from obtaining for $\mathcal{I}$ 
the same H\"older regularity as the control $X$. Let us make this rigorous, as follows. 

\vspace{0.3cm}

Let us first consider the linear case, that is we want to define the Young integral 
\beq
\mathcal{I}(s,t)=\int_0^s\int_\re G_{s-u}(t,v) X(du,dv).
\label{5.1}
\eeq
We have the following result on existence and regularity of the integral $\mathcal{I}$. 
\begin{proposition}
Recall that $\mathcal{C}^{\gam,\hat{\gam}}$ has been defined by
relation (\ref{defcgag}), and
assume that $X$ belongs to the space $\mathcal{C}^{\gam,\hat{\gam}}$, with $\gam+\hat{\gam}>1$. Then, the expression (\ref{5.1}) is well 
defined in the Young sense. Moreover, the function $\{\mathcal{I}(s,t), (s,t)\in [0,T]\times \re\}$ belongs to $\mathcal{C}^{\eta,\hat{\eta}}$, with $\eta=\rho\hat{\gam}+\gam-1$ and $\hat{\eta}=(1-\rho)\hat{\gam}$, for all $\rho\in \left(\frac{1-\gam}{\hat{\gam}},1\right)$.
\label{proposition2}
\end{proposition} 

\begin{remark} 
According to the preceding result, we deduce that $\eta+\hat{\eta}=\gam+\hat{\gam}-1$. Thus, the integral $\mathcal{I}$ has substantially lost regularity with respect to the control $X$. Furthermore, here and in the remainder of
the section, the Young integrals we consider are not really defined
as such, since they are based on a particular partition, suitable for
computations. We have chosen that solution for sake of simplicity, but
we believe our calculations could be carried out for a general 
family of partitions too.
\label{remark3}
\end{remark} 

\vspace{0.3cm}

Let us sketch the proof of the above proposition.

\vspace{0.3cm}

\noindent
{\it{Proof of Proposition \ref{proposition2}}}. We use a Young approximation to define the integral in the right-hand side of (\ref{5.1}). 
Fix $(s,t)\in [0,T]\times \re$ and consider the following dyadic grid on the rectangle in $\re^2$ determined by 
$\{(0,t-s),(s,t-s),(0,t+s),(s,t+s)\}$:
$$u^n_i=s\frac{i}{2^n}, \quad i=0,1,\dots,2^n,$$
$$v^n_j=t-s+s\frac{j}{2^n}, \quad j=0,1,\dots,2^{n+1}.$$
Notice that the above rectangle contains the support of the function 
$(s,y)\mapsto G_{s-u}(t,v)$ and each of the rectangles of the grid is 
of size $\frac{s}{2^n}\times \frac{s}{2^n}$. We set also
$$I^n_{i,j}=[u^n_i,u^n_{i+1}]\times [v^n_j,v^n_{j+1}]$$
and denote by $G_{s,t}(u,v)$ the function $G_{s-u}(t,v)\1_{[0,s]}(u)$. 

\vspace{0.3cm}

The natural Riemann type sum approximating the integral in (\ref{5.1}) would be
$$J_n(s,t)=\sum_{i=0}^{2^n-1} \sum_{j=0}^{2^{n+1}-1} G_{s,t}(u^n_i,v^n_j)\del_{I^n_{i,j}}X.$$
Equivalently, we study the convergence of the series
\beq
\sum_{n=1}^\infty \left(J_{n+1}(s,t)-J_n(s,t)\right),
\label{5.2}
\eeq
which we will decompose in the following way: first, notice that
$$\Del_{I^n_{i,j}}X = \Del_{I^{n+1}_{2i,2j}}X + \Del_{I^{n+1}_{2i+1,2j}}X + \Del_{I^{n+1}_{2i,2j+1}}X + 
\Del_{I^{n+1}_{2i+1,2j+1}}X.$$
On the other hand, the term $J_{n+1}(s,t)$ can be written as
\begin{align*}
J_{n+1}(s,t) =& \sum_{i=0}^{2^n-1} \sum_{j=0}^{2^{n+1}-1} \left[ G_{s,t}(u^{n+1}_{2i},v^{n+1}_{2j})\Del_{I^{n+1}_{2i,2j}}X + 
G_{s,t}(u^{n+1}_{2i+1},v^{n+1}_{2j})\Del_{I^{n+1}_{2i+1,2j}}X \right.\\
&\quad + \left. G_{s,t}(u^{n+1}_{2i},v^{n+1}_{2j+1})\Del_{I^{n+1}_{2i,2j+1}}X + 
G_{s,t}(u^{n+1}_{2i+1},v^{n+1}_{2j+1})\Del_{I^{n+1}_{2i+1,2j+1}}X\right].
\end{align*}
Thus
$$
J_{n+1}(s,t)-J_n(s,t) = A_1(s,t,n)+A_2(s,t,n)+A_3(s,t,n),
$$
where
\begin{align*}
A_1(s,t,n)= & \sum_{i=0}^{2^n-1} \sum_{j=0}^{2^{n+1}-1} \left( G_{s,t}(u^{n+1}_{2i},v^{n+1}_{2j})- 
G_{s,t}(u^{n+1}_{2i+1},v^{n+1}_{2j})\right) \Del_{I^{n+1}_{2i+1,2j}}X,\\
A_2(s,t,n)= & \sum_{i=0}^{2^n-1}  \sum_{j=0}^{2^{n+1}-1} \left( G_{s,t}(u^{n+1}_{2i},v^{n+1}_{2j})- 
G_{s,t}(u^{n+1}_{2i},v^{n+1}_{2j+1})\right) \Del_{I^{n+1}_{2i,2j+1}}X,\\
A_3(s,t,n)= & \sum_{i=0}^{2^n-1}  \sum_{j=0}^{2^{n+1}-1} \left( G_{s,t}(u^{n+1}_{2i},v^{n+1}_{2j})- 
G_{s,t}(u^{n+1}_{2i+1},v^{n+1}_{2j+1})\right) \Del_{I^{n+1}_{2i+1,2j+1}}X.
\end{align*}
Next, using the fact that $G_{s,t}$ is the indicator function of the light 
cone, a possible strategy for the estimation of those terms
is the following:
\begin{enumerate}
\item
Identify the non-trivial contributions 
in the sums $A_k(s,t,n)$, $k=1,2,3$, and then bound them in a suitable way.
\item
In order to get the desired bounds for the non-trivial contributions,
use the fact that
$$
|\Del_R X| \leq C \left(\frac{1}{2}\right)^{m(\gam+\hat{\gam})}
$$
for a rectangle $R$ of order $m$ in the grid defined above.
\end{enumerate}
Using these ideas, one can obtain the following estimation:
$$
|J_{n+1}(s,t)-J_n(s,t)|\leq  \frac{C}{2^{n(\gam+\hat{\gam}-1)}},
$$
with some positive constant $C$ independent of $n$. Since by hypothesis $\gam+\hat{\gam}>1$, 
we deduce that the series (\ref{5.2}) is convergent, and its sum
defines properly the integral in the right-hand side of (\ref{5.1}).

\vspace{0.3cm}

Let us outline the main computations leading to the regularity result
for the function $\{\mathcal{I}(s,t), (s,t)\in [0,T]\times \re\}$:
we fix $0\leq s_1<s_2\leq T$ and $-\infty<t_1<t_2<+\infty$, and we will 
study the quantity 
$|\Del_R \mathcal{I}|$, where $R=[s_1,s_2]\times [t_1,t_2]$. 
Furthermore, 
we assume that $t_1+s_2<t_2-s_2$, which implies that the light cones with vertices $(s_2,t_1)$ and $(s_2,t_2)$ do not intersect (the other case can be treated analogously). 
With  these assumptions in mind, we have
$$
\Del_R \mathcal{I}=
\mathcal{I}(s_2,t_2)-\mathcal{I}(s_2,t_1)-\mathcal{I}(s_1,t_2)
+\mathcal{I}(s_1,t_1)
$$
or, equivalently, 
\begin{equation}\label{idifg}
\int_0^T\int_\re \left( G_{s_2,t_2}(u,v)-G_{s_2,t_1}(u,v)-G_{s_1,t_2}(u,v)
+G_{s_1,t_1}(u,v)\right) X(ds,dy).
\end{equation}
We set 
$$D(u,v)=G_{s_2,t_2}(u,v)-G_{s_2,t_1}(u,v)-G_{s_1,t_2}(u,v)+G_{s_1,t_1}(u,v),
$$
and we proceed to the analysis of (\ref{idifg}) along the same lines
as in the first part of the proof. Namely, we will use Young type 
approximations of the following form:
fix $n_0\in \mathbb{N}$, and let $i_0\in \mathbb{Z}_-$ be such that 
$$u^{n_0}_{i_0}:=  s_1+(s_2-s_1)\frac{i_0}{2^{n_0}}<0\quad\textrm{and}\quad u^{n_0}_{i_0+1}:=s_1+(s_2-s_1)\frac{i_0+1}{2^{n_0}}\geq 0.$$
We consider the dyadic grid on $[0,T]\times \re$ given by $\Pi^n=\Pi^n_1\uplus\Pi^n_2$, where   
$\Pi^n_l$ is the grid on the rectangle in $\re^2$ determined by $\{(u^{n_0}_{i_0},t_l-s_2+u^{n_0}_{i_0}),(s_2,t_l-s_2+u^{n_0}_{i_0}),(u^{n_0}_{i_0},t_l+s_2-u^{n_0}_{i_0}),(s_2,t_l+s_2-u^{n_0}_{i_0})\}$, $l=1,2$, in the 
following way:
\begin{align*}
u^n_{i}=& \; s_1+(s_2-s_1)\frac{i}{2^n},\quad i=2^{n-n_0}i_0,\dots,2^n,\\
v^{l,n}_j= & \; t_l+(s_2-s_1)\frac{j}{2^n},\quad j=-2^n(1-2^{-n_0}i_0),\dots,0,\dots,2^n(1-2^{-n_0}i_0),\\
\end{align*}
for all $n\geq n_0$. Observe that
each square of the grid $\Pi^n_l$ is of size $\frac{s_2-s_1}{2^n}\times \frac{s_2-s_1}{2^n}$, $l=1,2$. Set now
$$I^{l,n}_{i,j}=[u^n_i,u^n_{i+1}]\times [v^{l,n}_j,v^{l,n}_{j+1}],$$
and notice that, once again, we will try to get some information
about the convergence of the series
\begin{equation}\label{serj}
\sum_{n=n_0}^{\infty}(J_{n+1}(s_1,s_2;t_1,t_2)-J_n(s_1,s_2;t_1,t_2)),
\end{equation}
where
$$J_n(s_1,s_2;t_1,t_2)=\sum_{\stackrel{i,j,l}{(u^n_i,v^{l,n}_j)\in\Pi^n}} D(u^n_i,v^{l,n}_j)\Del_{I^{l,n}_{i,j}}X.$$
Now, the desired convergence can be obtained in the same
spirit as for the first part of the proof, but with slightly more effort. 
For sake of conciseness, the details are left to the reader, but
let us briefly justify the need of introducing the parameter $\rho$. Namely, it turns out that one of the terms that we have to deal with 
for the convergence of (\ref{serj}) is given by
$$A(n)=\sum_{k=-i_n}^{2^n-i_n-1}\left( \Del_{I^{2,n+1}_{2i_n+1,2k}}X - \Del_{I^{1,n+1}_{2i_n+1,2k}}X\right),$$
where $i_n$ is some integer such that $0$ is contained in $(u^n_{i_n},u^{n+1}_{2i_n+1}]$. In order to get bounds on $A(n)$ of the form 
$(s_2-s_1)^{\mu}(t_2-t_1)^{\hat{\mu}}$, $\mu, \hat{\mu}\in (0,1)$, we 
notice that 
$$\Del_{I^{2,n+1}_{2i_n+1,2k}}X - \Del_{I^{1,n+1}_{2i_n+1,2k}}X = 
\Del_{R^n_k}X - 
\Del_{R^n_{k+1/2}}X,$$
where
$$R^n_k=\left[u^{n+1}_{2i_n+1},u^{n+1}_{2i_n+2}\right]\times \left[y^{1,n+1}_{2k},y^{2,n+1}_{2k}\right].$$
Then, owing to the hypothesis on $X$ and the definition of $(u^n_i,v^{l,n}_j)$, for any $\rho \in (0,1)$, we obtain
\begin{align*}
|A(n)|\leq  & \sum_{k=-i_n}^{2^n-i_n-1}\left| \Del_{I^{2,n+1}_{2i_n+1,2k}}X - \Del_{I^{1,n+1}_{2i_n+1,2k}}X\right|^\rho \left|\Del_{R^n_k}X - 
\Del_{R^n_{k+1/2}}X\right|^{1-\rho} \\
& \leq C 2^n \left( \frac{s_2-s_1}{2^{n+1}}\right)^{\rho(\gam+\hat{\gam})} 
\left( \frac{s_2-s_1}{2^{n+1}}\right)^{(1-\rho)\gam} (t_2-t_1)^{(1-\rho)\hat{\gam}} \\
&\leq C \left( \frac{1}{2}\right)^{n(\rho\hat{\gam}+\gam -1)}
(s_2-s_1)^{\rho\hat{\gam}+\gam}(t_2-t_1)^{(1-\rho)\hat{\gam}}.  
\end{align*}
This bound is then used to obtain the desired regularity in our 
statement. \qed

\vspace{0.5cm}

Let us jump to the non-linear case, which amounts to give a rigorous 
meaning and obtain regularity properties of the integral
$$\mathcal{I}_Z(s,t)=\int_0^s \int_{\re} G_{s-u}(t,v) Z(u,v)X(du,dv),$$
for a given process $Z$. More specifically,
we consider the following hypothesis for the function $Z$:

\vspace{0.3cm}

\noindent $\mathbf{ (H_\al)}$ 
The function
$Z:[0,T]\times \re \rightarrow \re$ is bounded, $Z(0,t)=0$, for all $t\in \re$, and there exist $\al\in (0,1)$ and a positive constant $C$ such that
$$|Z(s_1,t)-Z(s_2,t)| \leq C |s_1-s_2|^\al,$$
for any $s_1,s_2\in[0,T]$, $t\in \re$.

\vspace{0.3cm}

\noindent
Then we get the following existence and regularity result:
\begin{proposition}
Assume that $X$ belongs to the space $\mathcal{C}^{\gam,\hat{\gam}}$, with $\gam+\hat{\gam}>\frac{5}{3}$. 
We suppose also that the function $Z$ belongs to $\mathcal{C}^{\rho\hat{\gam}+\gam-1,(1-\rho)\hat{\gam}}$, with $\rho\in \left(\frac{3-2\gam-\hat{\gam}}{\hat{\gam}},\frac{2\hat{\gam}+\gam-2}{\hat{\gam}}\right)$, and that hypothesis $(H_\al)$ is satisfied for all $\al\in (2(2-\gam-\hat{\gam}),1)$. Then, the integral $\mathcal{I}_Z (s,t)$ is well defined in the Young sense and it defines a function belonging 
to $\mathcal{C}^{\rho\hat{\gam}+\gam-1,(1-\rho)\hat{\gam}}$.
\label{proposition4}
\end{proposition}

\begin{remark}
Here again, for sake of conciseness, we will omit the proof of this 
proposition, which is much more involved than the one of Proposition 
\ref{proposition2}, but may be carried out using analogous arguments.
\end{remark}

\begin{remark}
The motivation to construct the pathwise integral $\mathcal{I}_Z$ is to solve, given some control $X$, the 
integral  equation 
$$Y(s,t)=\int_0^s \int_{\re} G_{s-u}(t,v) \sig(Y(u,v))X(du,dv).$$
To prove existence and uniqueness of solution to the above equation, one usually applies a fixed point argument or, equivalently, a Picard iteration scheme. For this reason, in the preceding Proposition \ref{proposition4} we have assumed that the function $Z$ has the same regularity as the integral in the linear case (see Proposition \ref{proposition2}).
\end{remark}

\begin{remark}
According to the statement of Proposition \ref{proposition4}, we deduce that not only the integral $\mathcal{I}_Z$ looses again regularity with respect to the control $X$, but also the regularity of the latter must be strengthened,
since we assume now  $\gamma+\hat\gamma>5/3$ instead of $\gamma+\hat\gamma>1$.
\end{remark}


\end{document}